\documentclass[11pt]{article}

	\usepackage{amsmath, amssymb, amscd, latexsym, epsfig}
\newcommand{\arrow}{\rightarrow}
\newcommand{\cD}{{\cal D}}

\newcommand{\bb}{\mathbb}
\newcommand{\cx}{{\bb C}}

\newcommand{\half}{{\bb H}}
\newcommand{\integers}{{\bb Z}}

\newcommand{\reals}{{\bb R}}

\newcommand{\makefig}[3]{
	\begin{figure}[htbp]
        \refstepcounter{figure}
	\label{#2}
        \begin{center}
		~#3~\\
		\medskip
                {\sf Figure \thefigure.  #1}
        \end{center}
	\end{figure}
}

\thicklines
\newcommand{\maketab}[3]{
	\begin{figure}[htbp]
        \refstepcounter{figure}
	\label{#2}
        \begin{center}
		#3~\\
		\bigskip
                {\sf Table \thefigure.  #1}
        \end{center}
	\end{figure}
}

\renewcommand{\bold}[1]{\smallskip \noindent {\bf \boldmath #1 }\nopagebreak[4]}

\newcommand{\qed}{\nopagebreak[4]\hfill
\rule{2mm}{2.5mm} \bigskip \pagebreak[2]}

\renewcommand{\tilde}{\widetilde}

\newcommand{\asyto}{\sim}

\newcommand{\bs}{\backslash}

\newcommand{\bijects}{\leftrightarrow}

\newcommand{\brackets}[1]{\langle #1 \rangle}

\newcommand{\congruent}{\equiv}

\newcommand{\del}{\partial}

\newcommand{\equi}{\sim}

\newcommand{\isom}{\cong}

\newcommand{\mem}{\in}

\newcommand{\plusorminus}{\pm}

\newcommand{\AND}{\;\;\;\text{and}\;\;\;}

\newcommand{\st}{\: : \:}         

\newcommand{\abar}{{\overline{a}}}
\newcommand{\bbar}{{\overline{b}}}
\newcommand{\cbar}{{\overline{c}}}

\newcommand{\gbar}{{\overline{g}}}

\newcommand{\delbar}{{\overline{\del}}}

\newcommand{\chat}{{\widehat{\cx}}}

\newcommand{\Aut}{\operatorname{Aut}}

\newcommand{\id}{\operatorname{id}}

\renewcommand{\Im}{\operatorname{Im}}

\newcommand{\Isom}{\operatorname{Isom}}

\newcommand{\Ker}{\operatorname{Ker}}

\renewcommand{\mod}{\operatorname{mod}}

\newcommand{\SL}{\operatorname{SL}}

\newcommand{\zed}{\integers}

\newcommand{\cA}{{\cal A}}

\newcommand{\cG}{{\cal G}}
\newcommand{\cM}{{\cal M}}

\newtheorem{theorem}{Theorem}[section]
\newtheorem{prop}[theorem]{Proposition}
\newtheorem{lemma}[theorem]{Lemma}
\newtheorem{cor}[theorem]{Corollary}

\makeatletter
\def\cleardoublepage{\clearpage\if@twoside \ifodd\c@page\else
    \thispagestyle{plain}\hbox{}\newpage\if@twocolumn\hbox{}\newpage\fi\fi\fi}

\def\ps@headings{\let\@mkboth\markboth
  \def\@oddfoot{}%
  \def\@evenfoot{}%
  \def\@evenhead{\small \sc\thepage\hfil\leftmark}
  \def\@oddhead{\small \sc \rightmark\hfil\thepage}
  \def\chaptermark##1{{
    \edef\@tempa{\ifnum \c@secnumdepth >\m@ne \@chapapp\ \thechapter. \fi}%
    \expandafter \markboth \expandafter{\@tempa ##1}{}}}%
  \def\schaptermark##1{\markboth {##1}{##1}}%
  \def\sectionmark##1{{
    \edef\@tempa{\ifnum \c@secnumdepth >\z@ \thesection. \fi}%
    \expandafter \markright \expandafter{\@tempa ##1}}}}

\def\thebibliography#1{\section*{References\@mkboth
 {References}{References}}\list
 {[\arabic{enumi}]}{\settowidth\labelwidth{[#1]}\leftmargin\labelwidth
 \advance\leftmargin\labelsep
 \usecounter{enumi}}
 \def\newblock{\hskip .11em plus .33em minus .07em}
 \sloppy\clubpenalty4000\widowpenalty4000
 \sfcode`\.=1000\relax}

\newif\if@restonecol
\def\theindex{\@restonecoltrue\if@twocolumn\@restonecolfalse\fi
\columnseprule \z@
\columnsep 35pt\twocolumn[\@makeschapterhead{Index}]
 \@mkboth{Index}{Index}\thispagestyle{plain}\parindent\z@
 \parskip\z@ plus .3pt\relax\let\item\@idxitem}
\def\@idxitem{\par\hangindent 40pt}

\def\endtheindex{\if@restonecol\onecolumn\else\clearpage\fi}

\def\footnoterule{\kern-3\p@ 
 \hrule width .4\columnwidth 
 \kern 2.6\p@} 
\@addtoreset{footnote}{chapter} 
\long\def\@makefntext#1{\parindent 1em\noindent 
 \hbox to 1.8em{\hss$^{\@thefnmark}$}#1}

\@addtoreset{equation}{section}

\renewcommand{\l@section}{\@dottedtocline{0}{1.5em}{2.3em}}
\renewcommand{\l@subsection}{\@dottedtocline{1}{3.8em}{3.2em}}
\renewcommand{\l@subsubsection}{\@dottedtocline{2}{7.0em}{4.1em}}

\makeatother

\newcommand{\Gt}{\tilde{G}}
\newcommand{\Mt}{\tilde{M}}
\newcommand{\gammat}{\tilde{\gamma}}

\begin{document}

\title{ \vspace{-1in}
        {\bf Almost simple geodesics on the triply punctured sphere}
\vspace{.2in}}

\author{Moira Chas, Curtis T. McMullen and Anthony Phillips }

\date{5 March 2017}

\maketitle

\begin{abstract}
Every closed hyperbolic geodesic
$\gamma$ on the triply--punctured sphere 
$M = \chat - \{0,1,\infty\}$ has a
self--intersection number $I(\gamma) \ge 1$ and a
combinatorial length $L(\gamma) \ge 2$, the latter defined by
the number of times $\gamma$ passes through the upper
halfplane.

In this paper we show that $\delta(\gamma) = I(\gamma) - L(\gamma) \ge -1$
for all closed geodesics; and that for each fixed $\delta$,
the number of geodesics with invariants $(\delta,L)$ 
is given exactly by a quadratic polynomial $p_\delta(L)$ for all $L \ge 4 + \delta$.
\end{abstract}

\tableofcontents

\vfill \footnoterule \smallskip
{\footnotesize \noindent
        Research supported in part by the NSF.
}

\thispagestyle{empty}
\setcounter{page}{0}

\newpage

\section{Introduction}
\label{sec:intro}

There are no simple closed geodesics on the triply--punctured sphere.
That is, the geometric self--intersection number 
$I(\gamma)$ of every
closed hyperbolic geodesic $\gamma$ on the Riemann surface 
\begin{displaymath}
	M = \chat - \{0,1,\infty\} 
\end{displaymath}
(endowed with its complete conformal metric of constant curvature $-1$)
satisfies $I(\gamma) > 0$.

In the absence of simple loops, one can aim instead
to classify and enumerate those
geodesics on $M$ that are {\em almost simple}, in the sense that $I(\gamma)$ is
small compared to the {\em combinatorial length} $L(\gamma)$.
For our purposes, it is convenient to define
$L(\gamma)$ to be the number of times that
$\gamma$ passes through the upper halfplane;
equivalently, $2L(\gamma)$ is the number of times that $\gamma$ crosses the real line $\reals \subset \chat$ (for an example, see Figure \ref{fig:geod_example}).
Our first result (\S\ref{sec:minusone}) relates these two quantities.

\begin{theorem}
\label{thm:minusone}
For any closed geodesic $\gamma \subset M$, we have $I(\gamma) \ge L(\gamma) - 1$.
\end{theorem}

The {\em defect} 
\begin{displaymath}
	\delta(\gamma) = I(\gamma) - L(\gamma) \ge -1
\end{displaymath}
is thus a natural measurement of the failure of $\gamma$ to be simple.  For a typical long geodesic, $\delta(\gamma)$ is on the order
of $L(\gamma)^2$.  We will be interested in the opposite regime,
where $\delta(\gamma) = O(1)$.
More precisely, for each fixed $\delta$ we wish to study the function
\begin{displaymath}
	N_\delta(L) = |\{\gamma \subset M \st L(\gamma) = L \AND \delta(\gamma) = \delta \} | .
\end{displaymath}
Here we have identified $\gamma$ with a subset of $M$,
so $N_\delta(L)$ is a count of the number of 
{\em unoriented, primitive} geodesics of length $L$ and defect $\delta$.
Table \ref{tab:main} gives the value of $N_\delta(L)$ for small $\delta$ and $L$.

\makefig{A closed hyperbolic geodesic with $L(\gamma)=4$ and $I(\gamma)=5$.}{fig:geod_example}{
\includegraphics[height=2.1in]{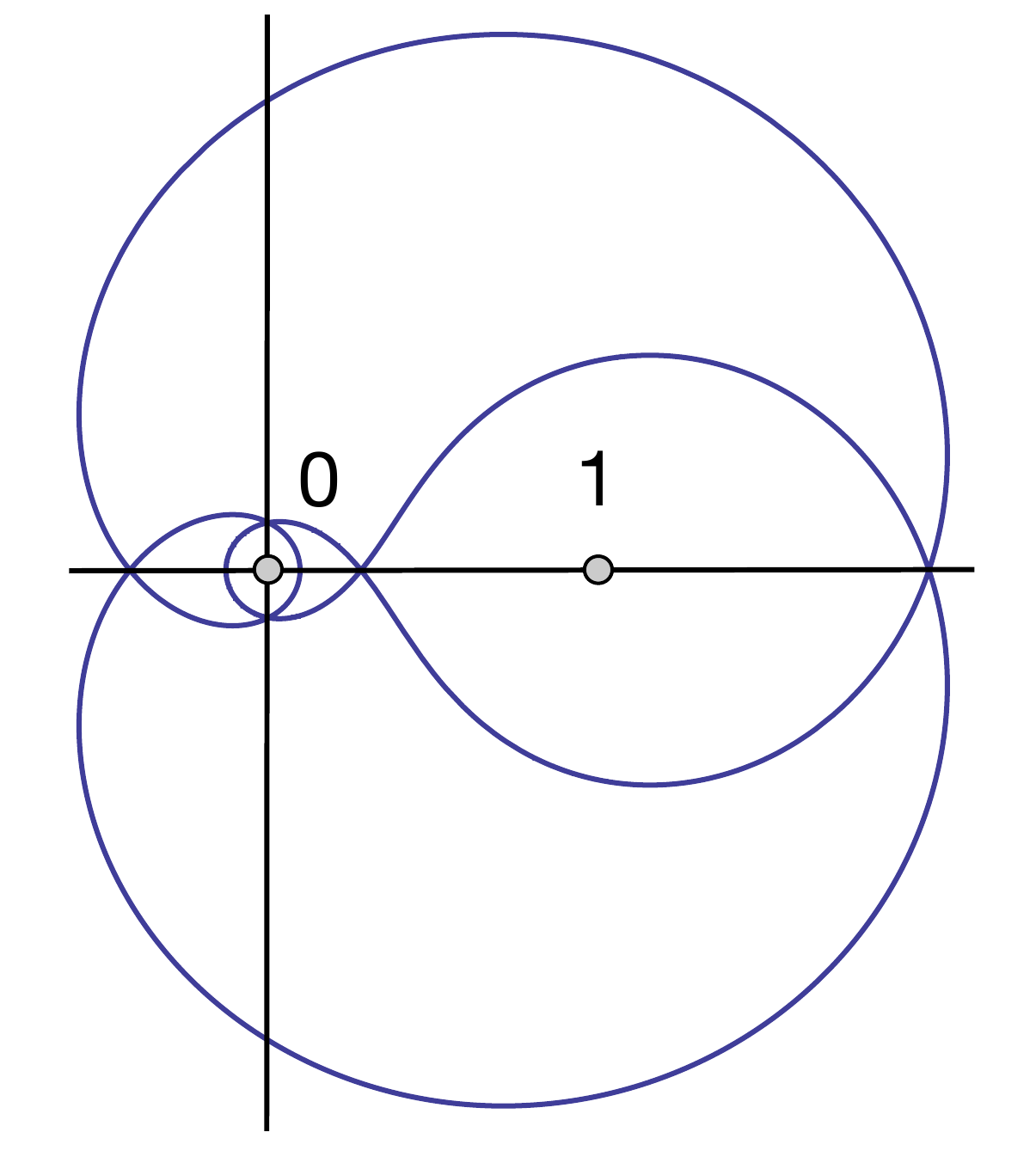}}

In this paper we show:

\begin{theorem}[Quadratic enumeration]
\label{thm:main}
For each $\delta \ge -1$, there exists a quadratic polynomial $p_\delta(L)$ such that
$N_\delta(L) = p_\delta(L)$ for all $L \ge \delta+4$.
\end{theorem}

We emphasize that Theorem \ref{thm:main} concerns
the exact value of $N_\delta(L)$, not just its asymptotic behavior.
The polynomials $p_\delta(L)$ for $\delta \le 11$ can be found
by examining the columns of Table \ref{tab:main};
for example, we have:
\begin{eqnarray*}
	p_{-1}(L) & = & 3L^2 - 9L + 9, \\
	p_{0}(L) & = & 4L^2 - 24L + 38, \\
	p_{1}(L) & = & 30L^2 - 240 L + 486,
	\AND \\
	p_{11}(L) & = & 16608L^2 - 363900L + 2030832 . 
\end{eqnarray*}
The statements of Theorems \ref{thm:minusone} and
\ref{thm:main} were first suggested by the experimental data
in this table.

\bold{Question.}
Is it true that $N_\delta(L)=0$  if and only if $\delta > (L-3)(L-1)/3$?

\bold{Pairs of pants.}
Although stated in terms of geodesics,
Theorems \ref{thm:minusone} and \ref{thm:main}
can be regarded as topological results about closed loops 
on $M$ (or equivalently, on a pair of pants).
Indeed, every essential, nonperipheral closed loop on $M$ 
is homotopic to a unique geodesic $\gamma$,
and $I(\gamma)$ is simply the minimum number 
of (transverse) self-intersections among all representatives of that 
homotopy class.
The combinatorial length $L(\gamma)$ can also be described topologically,
in terms of generators for $\pi_1(M)$ (see \S\ref{sec:background}).  
For more details on the geometric intersection number,
see e.g. \cite{FLP}, \cite{Bonahon:currents}.

\maketab{
	The number $N_\delta(L)$ of unoriented, primitive geodesics on $\chat - \{0,1,\infty\}$ 
	of combinatorial length $L$ and self-intersection number $I = L+\delta$.
	}{tab:main}{
\hspace{-1.5in}{
{\footnotesize 
\setlength{\tabcolsep}{0.045in}

\begin{tabular}{|r||rrrrrrrrrrrrr|}
	\cline{2-14}
\multicolumn{1}{c|}{}
    & $\delta = \,-1$ & 0 & 1 & 2 & 3 & 4 & 5 & 6 & 7 & 8 & 9 & 10 & 11 \\ 
	\cline{2-14}
	\cline{2-14}
	\hline
$L=2$ & 3 & 0 & 0 & 0 & 0 & 0 & 0 & 0 & 0 & 0 & 0 & 0 & 0 \\ \hline
3 & 9 & 1 & 0 & 0 & 0 & 0 & 0 & 0 & 0 & 0 & 0 & 0 & 0 \\ \hline
4 & 21 & 6 & 3 & 0 & 0 & 0 & 0 & 0 & 0 & 0 & 0 & 0 & 0 \\ \hline
5 & 39 & 18 & 36 & 9 & 0 & 0 & 0 & 0 & 0 & 0 & 0 & 0 & 0 \\ \hline
6 & 63 & 38 & 126 & 54 & 27 & 18 & 9 & 0 & 0 & 0 & 0 & 0 & 0 \\ \hline
7 & 93 & 66 & 276 & 156 & 216 & 150 & 135 & 51 & 21 & 6 & 0 & 0 & 0 \\ \hline
8 & 129 & 102 & 486 & 318 & 666 & 528 & 672 & 438 & 375 & 180 & 78 & 72 & 36 \\ \hline
9 & 171 & 146 & 756 & 540 & 1386 & 1218 & 2070 & 1648 & 1995 & 1269 & 1088 & 948 & 660 \\ \hline
10 & 219 & 198 & 1086 & 822 & 2376 & 2226 & 4560 & 4044 & 5970 & 4632 & 5532 & 4890 & 4596 \\ \hline
11 & 273 & 258 & 1476 & 1164 & 3636 & 3552 & 8160 & 7764 & 13302 & 11571 & 16608 & 15342 & 18081 \\ \hline
12 & 333 & 326 & 1926 & 1566 & 5166 & 5196 & 12870 & 12818 & 24414 & 22806 & 36779 & 35838 & 49428 \\ \hline
13 & 399 & 402 & 2436 & 2028 & 6966 & 7158 & 18690 & 19206 & 39336 & 38574 & 67836 & 68925 & 105708 \\ \hline
14 & 471 & 486 & 3006 & 2550 & 9036 & 9438 & 25620 & 26928 & 58068 & 58890 & 110454 & 115806 & 191337 \\ \hline
15 & 549 & 578 & 3636 & 3132 & 11376 & 12036 & 33660 & 35984 & 80610 & 83754 & 164678 & 176844 & 309132 \\ \hline
16 & 633 & 678 & 4326 & 3774 & 13986 & 14952 & 42810 & 46374 & 106962 & 113166 & 230508 & 252060 & 460080 \\ \hline
17 & 723 & 786 & 5076 & 4476 & 16866 & 18186 & 53070 & 58098 & 137124 & 147126 & 307944 & 341454 & 644244 \\ \hline
\end{tabular}
}
	}
	}

\makefig{Topology of closed geodesics with $\delta(\gamma)=0$.}{fig:N0}{
\includegraphics[height=1.8in]{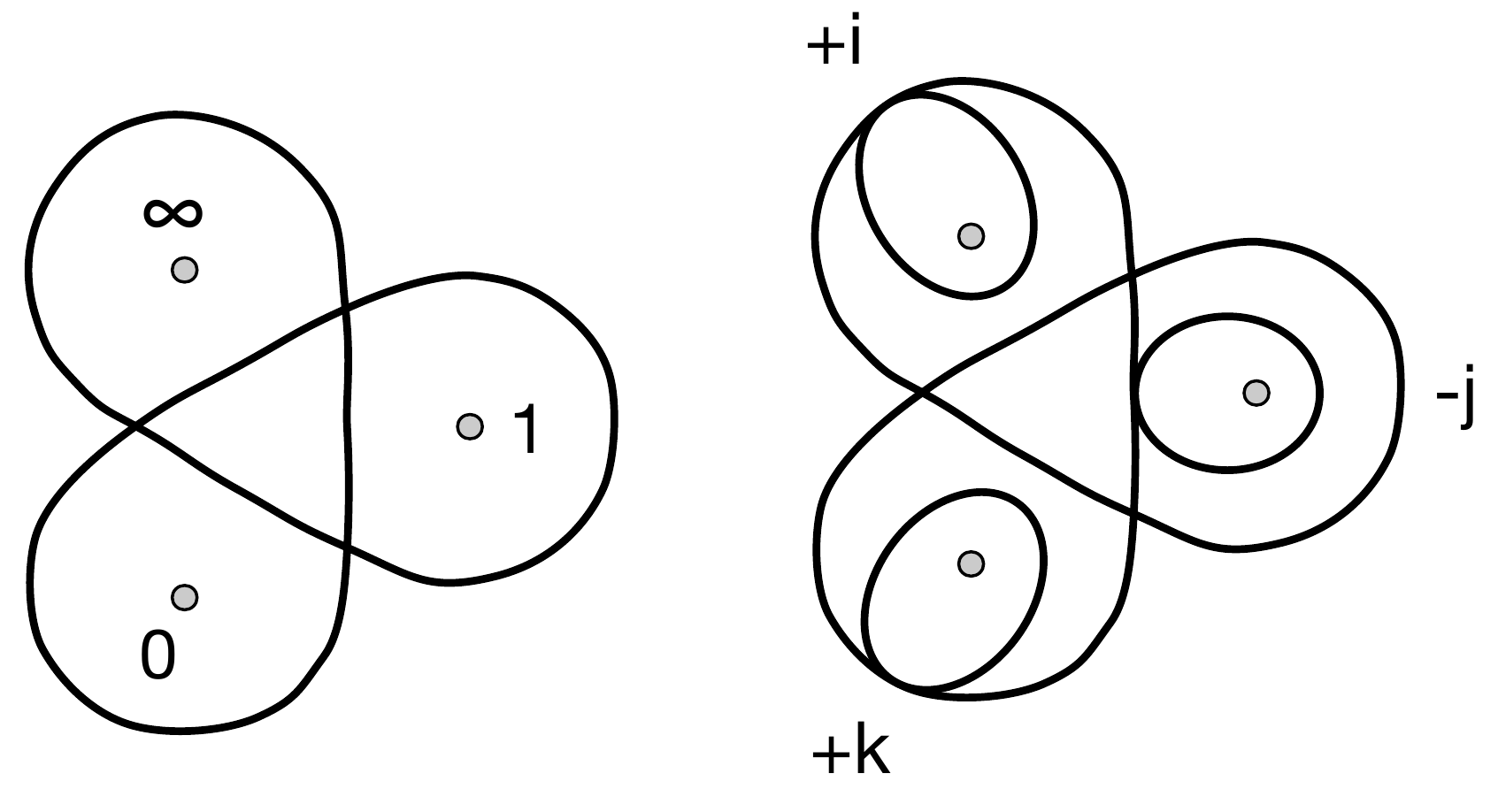}
}

\bold{Decorations.}
The mechanism 
behind Theorem \ref{thm:main} is illustrated in Figure \ref{fig:N0}.
The trefoil at the left in the figure shows a geodesic with $I(\gamma)=L(\gamma)=3$, and hence
$\delta(\gamma) = 0$.  For each $(i,j,k) \mem \zed^3$ we can decorate $\gamma$ 
by adding multiple loops around
each of the three punctures of $M$ to obtain the homotopy class of
another geodesic with $\delta(\gamma_{ijk}) = 0$; here the signs of 
$i$, $j$ and $k$ indicate if the decorations are to be attached to the inner or outer triangle of $\gamma$.
The homotopy class of $\gamma_{ijk}$ is indicated schematically
by the train--track shown at the right.
It turns out that every geodesic with $\delta = 0$ is obtained in this way, and hence $N_0(L)$
is simply the number of solutions to the equation
\begin{displaymath}
	L(\gamma_{ijk}) = 3 + |i|+|j|+|k|  = L.
\end{displaymath}
Explicitly, this count is given by
\begin{equation}
\label{eq:N0}
	N_0(L) = 8 \binom{L-4}{2} + 12 \binom{L-4}{1} + 6 \binom{L-4}{0}
		+ \binom{L-4}{-1} ,
\end{equation}
which agrees with the quadratic polynomial $p_0(L)$ for $L \ge 4$.

\bold{Binomial coefficients:  conventions.}
In the statement of this and other results, we adopt the convention that
\begin{equation}
\label{eq:bin1}
	\binom{n}{k} = 0 \; \text{if $n$ or $k$ is negative, except} \; \binom{-1}{-1} = 1 .
\end{equation}
This convention is chosen so that the equation
\begin{equation}
\label{eq:bin2}
	\binom{n+r-1}{r-1} = \left|\left\{(n_1,\ldots,n_r) \st n_i \ge 0, n_i \mem \zed, \sum_1^r n_i = n \right\}\right|
\end{equation}
is valid for all integers $n$ and all $r \ge 0$.
With this convention, the usual expression for $\binom{n}{k}$ as a 
polynomial in $n$ is valid only for $n \ge 0$.
For example, formula (\ref{eq:N0}) for $N_0(L)$ agrees with the 
quadratic polynomial $p_0(L)$ for $L \ge 4$, but for $L=1,2,3$ we have
$N_0(L) = 0,0,1$ while $p_0(L)=18,6,2$.

\bold{Motifs.}
More generally, to prove Theorem \ref{thm:main},
we will show that every closed geodesic in $M$ with $\delta(\gamma) = \delta$
is obtained by decorating one of finitely many {\em motifs}.
In terms of these motifs, we obtain a formula for $N_\delta(L)$ as a sum of binomial coefficients, valid for all $L$ (\S\ref{sec:motifs}):

\begin{theorem}[Binomial enumeration]
\label{thm:binom}
For all integers $\delta$ and $L$, we have
\begin{displaymath}
	N_\delta(L) = 
		\sum_{\text{motifs $\gamma$ with $\delta(\gamma)=\delta$}}
		\binom{L-L(\gamma) + \rho(\gamma)-1}{\rho(\gamma)-1} \cdot
\end{displaymath}
\end{theorem}
Here the {\em rank} of a motif, $0 \le \rho(\gamma) \le 3$, indicates how many decorations it admits.

\makefig{The closed geodesics $\gamma_{010}$, $\gamma_{101}$ and $\gamma_{1,-1,1}$.}{fig:3geods}{
\includegraphics[height=1.8in]{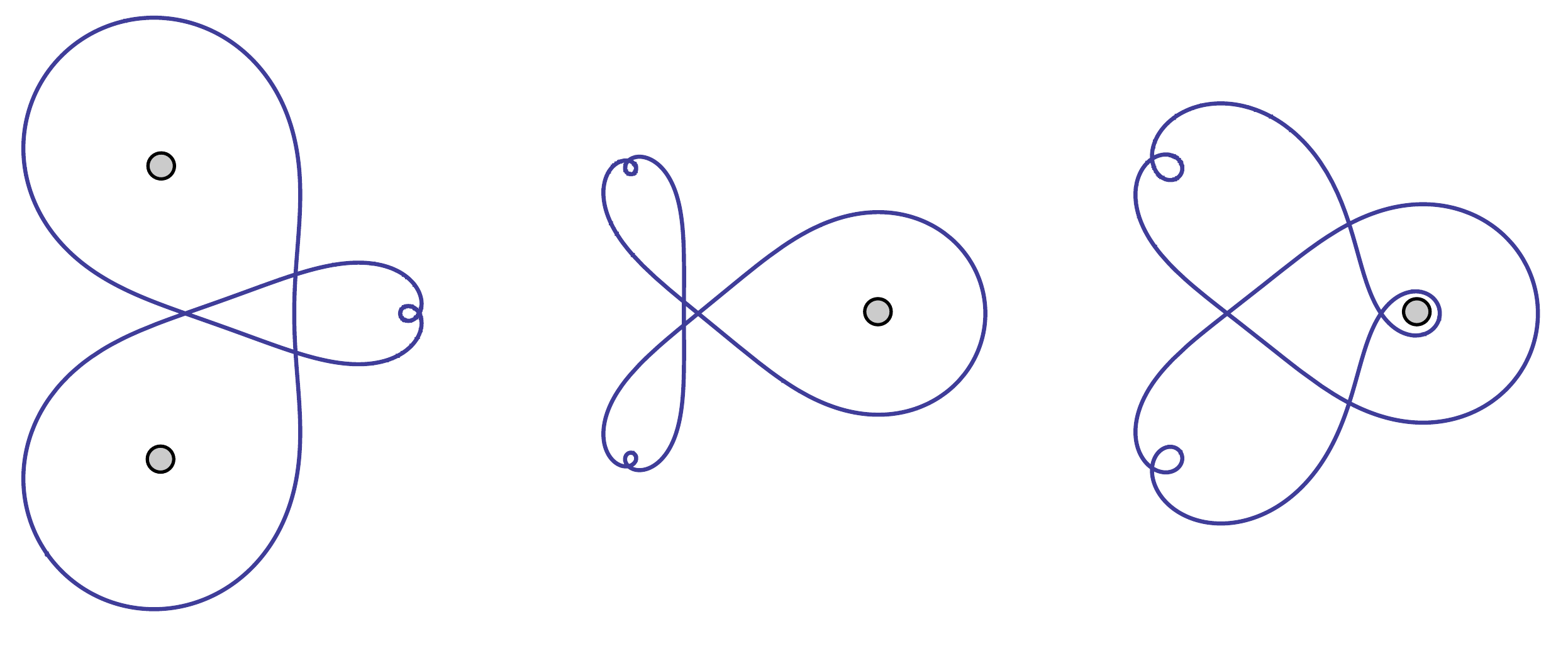}
}

For example, there are 27 motifs with $\delta(\gamma)=0$; they correspond, in terms of Figure \ref{fig:N0}, to the geodesics $\gamma_{ijk}$ 
with $i,j,k \mem \{-1,0,1\}$.
Three of these geodesics are shown in Figure \ref{fig:3geods}.
Grouping them together according their ranks, given by
$\rho(\gamma_{ijk}) = |i|+|j|+|k|$,
we obtain the 4 terms in equation
(\ref{eq:N0}) for $N_0(L)$.  
For more details on this example and others, see \S\ref{sec:motifs}.

\bold{Lengths of motifs.}
In \S\ref{sec:rank} we show
that the length of any motif satisfies the bound:
\begin{equation}
\label{eq:rankbound}
	\delta(\gamma) +\rho(\gamma) + 3 \ge L(\gamma).
\end{equation}
This crucial bound has two important consequences.  First, it implies
that Theorem \ref{thm:binom} expresses $N_\delta(L)$ as a {\em finite}
sum of binomial coefficients, since $6+\delta \ge L(\gamma)$;
and second, it shows that each binomial coefficient agrees with
a polynomial in $L$ of degree $\rho(\gamma) -1 \le 2$ in the range
$L \ge \delta+4$.  Thus $N_\delta(L)$ itself is a quadratic polynomial
in $L$, proving Theorem \ref{thm:main}.

We emphasize that Theorem \ref{thm:binom} gives a formula
for $N_\delta(L)$ that is valid for all $L$, not just
$L \ge \delta+4$.
Moreover, motifs have a simple description in terms of combinatorial
group theory (see \S\ref{sec:motifs}).  
The main subtlety in evaluating this formula comes from the 
condition $\delta(\gamma)=\delta$, which requires the computation of
the self--intersection number of $\gamma$.

\bold{Question.}
What is the behavior of the leading coefficient of $p_\delta(L)$?
This coefficient is one--half the number of motifs with defect $\delta$ and full rank.

\bold{Spheres with 4 or more punctures.}
Much of our analysis generalizes in a straightforward
way to the case where $M$ is an $n$--times punctured sphere, $n \ge 4$;
for example, Theorem \ref{thm:binom} remains valid in this setting.
The crucial difference is that for $n \ge 4$, there are infinitely many
motifs with $\delta(\gamma) = \delta$, so the formula for 
$N_\delta(L)$ becomes an infinite sum of binomial coefficients
(only finitely many of which are nonzero for a given value of $L$).

\bold{Perspectives on $M$.}
The proof of Theorem \ref{thm:main} pivots on two
crucial shifts in perspective on the geometry and combinatorial group theory
of the triply--punctured sphere $M$.

The first shift is to express words in $G = \pi_1(M,p)$ in terms
of generators, not for $G$, but for the reflection group
\begin{displaymath}
	\Gt = \brackets{x,y,z \st x^2=y^2=z^2=\id} 
\end{displaymath}
which contains $G$ with index two.
The advantage of these generators is that they do not
privilege the upper or lower halfplane, and that they allow one to
give combinatorial meaning to a geodesic which takes a half--integral
number of turns around a cusp.

The second shift is to replace the standard hyperbolic 
metric on $M$ with a complete hyperbolic metric of infinite volume,
turning the convex core of $M$ into a symmetric pair of pants with 
long boundary components.
In this new metric, the location of self--intersections
of geodesics is changed in an advantageous way,
even though the total number of self--intersections remains
the same.  
More precisely, we obtain a direct relationship between 
the {\em depth} of a geodesic excursion into one of the ends of $M$, and 
the {\em length} of a run of alternating letters in the 
corresponding word $w \mem \Gt$.

\bold{Outline of the paper.}
The interplay of hyperbolic geometry and combinatorial group theory
just described is used in \S\ref{sec:surgery}
to bound the change in $I(\gamma)$ when a loop 
around a cusp is added or removed.
A quick proof of Theorem \ref{thm:minusone} follows in 
\S\ref{sec:minusone}.
The theory of motifs is discussed in \S\ref{sec:motifs},
and the proof of Theorem \ref{thm:main} is completed in \S\ref{sec:rank}.

In the Appendix we give a formula for $I(\gamma)$
in terms of combinatorial group theory, suitable for computing
the entries in Table \ref{tab:main}.
We also establish the lower bound $I(\gamma) \ge L(\gamma)^2/6$
for a class of loops with controlled excursions into the cusps of $M$.
Both results play an important role 
in \S\ref{sec:rank}, where we prove the inequality (\ref{eq:rankbound})

\bold{Notes and references.}
For additional background on curves, surfaces and intersection numbers,
see e.g. \cite{FLP}, \cite{Hass:Scott:int}, 
\cite{Bonahon:currents}, \cite{Cohen:Lustig:sccs}, \cite{Stillwell:book:surfaces},
\cite{Despre:Lazarus:sccs} and the references therein.
Variants of Table \ref{tab:main}, 
not restricted to primitive geodesics, appear in \cite{Chas:Phillips:torus} and \cite{Chas:Phillips:triply}.
The statistical distribution of
self--intersection numbers is studied in \cite{Chas:Lalley:stat}.

\section{Background and notation}
\label{sec:background}

In this section we make explicit the relationship between closed geodesics and combinatorial group theory.

\bold{Group theory.}
We will work in the reflection group
\begin{displaymath}
	\Gt = \brackets{x,y,z \st x^2 = y^2 = z^2 = \id },
\end{displaymath}
and in the index two subgroup $G$ generated by $(a,b,c) = (xy,yz,zx)$; it has the presentation
\begin{displaymath}
	G = \brackets{a,b,c \st abc=\id} .
\end{displaymath}

A word $w$ in the generators $x,y,z$ of $\Gt$ is {\em reduced} if consecutive letters
of $w$ are distinct.  
Every element $w \mem \Gt$ is represented by a unique reduced word,
whose length will be denoted by $\ell(w)$.
We say $w$ is {\em cyclically reduced} if its first and last letters are also distinct.
The shortest words in a given conjugacy class are cyclically reduced.

We have $w \mem G$ iff $\ell(w)$ is even.
We say $w \mem G$ is {\em primitive} unless $w = g^n$ for some $n>1$ and $g \mem G$.

We say two cyclically reduced words $w_1,w_2 \mem G$ are {\em equivalent},
if $w_1$ is conjugate to
$w_2^{\plusorminus 1}$.  In terms of words, if $w_1 = g_1 \cdots g_{2n}$, then
\begin{displaymath}
	w_2 = ( g_{2i+1} \cdots g_{2n} g_1 \cdots g_{2i})^{\plusorminus 1}
\end{displaymath}
for some $i$.
Let
\begin{displaymath}
	\cG = \left\{w \mem G \st \text{$w$ is primitive, cyclically reduced and $\ell(w) \ge 4$}\right\} / \equi. 
\end{displaymath}
Each element of $\cG$ is an equivalence class $[\omega]$ of cyclically reduced words.
The length condition insures that $w$ involves all three generators $x,y,z$.

\bold{Hyperbolic geodesics.}
We now return to the triply--punctured sphere
\begin{displaymath}
	M = \chat - \{0,1,\infty\} .
\end{displaymath}
The Riemann surface $M$ can be presented as the quotient
of the upper halfplane $\half = \{z \mem \cx \st \Im(z) > 0\}$
by the group $\Gamma(2) \subset \SL_2(\zed)$ consisting of matrices
with $A \congruent I \mod 2$, acting by M\"obius transformations.
It carries a unique complete, conformal
{\em hyperbolic metric} of constant curvature $-1$, inherited
from the metric $|dz|/\Im(z)$ on $\half$.  Geodesics on $M$
are covered by semicircles perpendicular to the boundary in $\half$.

Choosing a basepoint $p \mem M$ with $\Im(p)>0$, we have an isomorphism
\begin{displaymath}
	\phi : \pi_1(M,p) \isom G .
\end{displaymath}
To define this map unambiguously, we label the components $(-\infty,0)$, $(0,1)$ and $(1,\infty)$ of $\reals \cap M$
by $X$, $Y$ and $Z$ respectively.  
Given a smooth loop $\gamma : [0,1] \arrow M$ transverse to the real axis, with $\gamma(0)=\gamma(1) = p$, 
let $0<t_1 < \cdots < t_{2n} < 1$ denote the parameters such that $\gamma(t_i) \mem \reals$.
To each such crossing we associate one the generators of $\Gt$, namely $g_i = x$ if $\gamma(t_i) \mem X$,
$g_i = y$ if $\gamma(t_i) \mem Y$, and $g_i = z$ if $t_i \mem Z$.  We then define
\begin{displaymath}
	\phi([\gamma]) = g_1 \cdots g_{2n} .
\end{displaymath}
It is readily verified that $\phi$ depends only on the homotopy class of $\gamma$ and gives an isomorphism to $G$.

This isomorphism determines a natural bijection
\begin{equation}
\label{eq:G}
	\cG \bijects \left\{\text{closed hyperbolic geodesics $\gamma \subset M$} \right\} .
\end{equation}
This map sends a conjugacy class in $\pi_1(M,p)$ to its geodesic representative.  Conversely, if we  traverse a closed geodesic $\gamma \subset M$
starting at a point in the upper halfplane, and write down the corresponding generator of $G$ each time $\gamma$ crosses the real axis,
we obtain a cyclically reduced word $w \mem G$ representing the corresponding class $[w] \mem \cG$.
(The elements of $\pi_1(M,p)$ with $\ell(w) \le 2$ are excluded from $\cG$ because they are peripheral or trivial.)

\bold{Self--intersection numbers of words.}
Using the bijection (\ref{eq:G}), we can regard length and self--intersection number
as functions of cyclically reduced words in $G$, defined by
$L(w) = L(\gamma_w)$ and $I(w) = I(\gamma_w)$ .

Clearly $L(w) = \ell(w)/2$.  A combinatorial algorithm for computing $I(w)$ directly in terms
of the word $w$ is described in the Appendix.

\section{Surgery}
\label{sec:surgery}

In this section we use hyperbolic geometry to prove some
purely combinatorial statements about the behavior of 
intersection numbers under elementary modifications of
a group element $w \mem G$.  

\bold{Overview.}
The main results of this section are
Theorems \ref{thm:int1}, \ref{thm:int2} and \ref{thm:int3} below.
They are used as tools in all the sections that follow.

Although the statements of these theorems 
are combinatorial, their proofs are geometric.
To carry out them out, we first open up the
cusps of $M$ to obtain a hyperbolic pair of pants
$M(\lambda)$, whose convex core is bounded by three geodesics
(cuffs) of equal length.

Each cuff of $M(\lambda)$ induces 
an annular covering space $\Mt(\lambda) \arrow M(\lambda)$ (Figure \ref{fig:covers}).
In this covering space, the
winding number $N(\gamma)$ of a geodesic lifted from
$M(\lambda) \isom M$ can be related to runs of
alternating letters in the corresponding word $w \mem \pi_1(M)$
(see Figure \ref{fig:wind3}).
On the other hand, the larger the winding number,
the closer $\gamma$ comes to one of the cuffs of $M(\lambda)$
(see equation (\ref{eq:Nd}) and Figure \ref{fig:annulus}). 
The part of a geodesic closest to the cuff can be decorated without adding
any unexpected new self--intersections, and the desired
bounds on $I(w)$ follow.

\bold{Runs.}
We turn to the statements of the main theorems of this section.

Consider a combinatorial closed geodesic $[w] \mem \cG$, represented by a cyclically reduced word $w$.
Since $w$ is primitive and $\ell(w) \ge 4$, all three generators $x,y,z$ occur as letters of $w$.

A {\em run} $r$ of $w$ is a maximal sequence of (cyclically) consecutive
letters of $w$ in which at most two generators appear.
(The location of $r$ in $w$ is part of the information determining a run.)
Clearly $2 \le \ell(r) < \ell(w)$.

If $r$ only involves the generators $x,y$, then we say $r$ is a
run of {\em type} $xy$.  Any run is of type $xy$, $yz$ or $zx$.

\bold{Expansion.}
Given a run $r$ of $w$, we can {\em expand} $r$ by repeating its first two letters to obtain
a word $r^+$ with $\ell(r^+) = \ell(r)+2$; and then expand $w$
by replacing $r$ with $r^+$.
The resulting word $w^+[r]$ is still cyclically reduced,
because the first and last letters of $r$ are the same
as those of $r^+$.  In particular, we have:
\begin{displaymath}
	L(w^+[r]) = L(w) + 1.
\end{displaymath}
For example, if $w=xyzy$, and $r$ is the initial run $xy$, then $w^+[r] = xyxyzy$.

\bold{Contraction.}
Similarly, given any run $r$ of $w$ with $\ell(r) \ge 3$, we can {\em contract} $r$ 
to $r^-$ by removing its first two letters, and then contract $w$ by replacing $r$ with $r^-$.
The resulting word $w^-[r]$ is still reduced, and hence
\begin{displaymath}
	L(w^-[r]) = L(w) - 1.
\end{displaymath}
For example, if $w = xzxyxy$, and $r$ is the initial run $xzx$, then
$w^-[r] = xyxy$.  This example shows that primitivity need not
be preserved by contraction.

\bold{Behavior of $I(w)$.}
The main results of this section control the behavior of the self--intersection number $I(w)$ under
expansion and contraction.  Let us say a run is {\em exceptional} if
\begin{displaymath}
	I(w^-[r]) = I(w) - 1.
\end{displaymath}

\begin{theorem}
\label{thm:int1}
A cyclically reduced word $w$ with $[w] \mem \cG$ has at most one exceptional run of each type.
\end{theorem}

\begin{theorem}
\label{thm:int2}
Let $r$ be a run of $w$ with $\ell(r) \ge 3$.  Then either:
\begin{displaymath}
	I(w^-[r]) = I(w) - 1,
\end{displaymath}
and $\ell(r) \ge \ell(s)$ for every other run $s$ of the same type as $r$; or
\begin{displaymath}
	I(w^-[r]) \le I(w)-3.
\end{displaymath}
\end{theorem}

\begin{theorem}
\label{thm:int3}
If $\ell(r) > \ell(s)$ for every other run $s$ of the same type as $r$, then $I(w^+[r]) = I(w)+1$.
\end{theorem}
Since contraction decreases the length of a run by two, we have:

\begin{cor}
If $\ell(r) > \ell(s)+2$ for every other run $s$ of the same type as $r$,
then $r$ is exceptional:  we have $I(w^-[r]) = I(w)-1$.
\end{cor}

\makefig{A hyperbolic geodesic in an annulus with winding number $N(\gamma)>1$.}{fig:annulus}{
\includegraphics[height=2.0in]{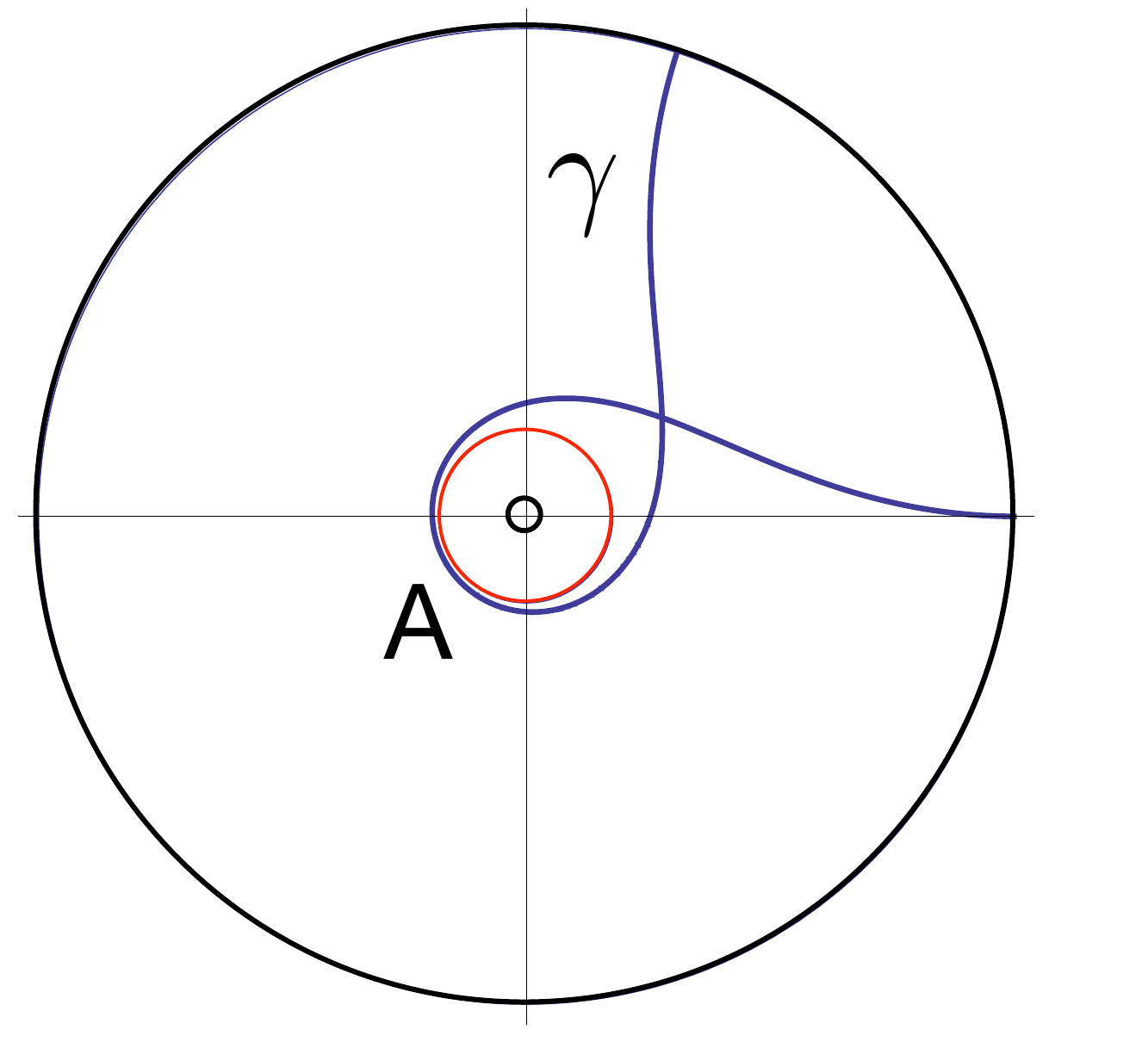}
}

\bold{Hyperbolic geodesics in annuli.}
The proofs of these results are based on hyperbolic geometry.
We begin by recalling some simple facts about hyperbolic geodesics in 
the annulus 
\begin{displaymath}
	U(r) = \{z \st r < |z| < 1 \}.
\end{displaymath}
Recall that the circle $A \subset U(r)$ defined by $|z| = \sqrt{r}$ is the unique
closed hyperbolic geodesic in the annulus.
If we denote the length of this geodesic by $\log(a)$, $a>1$, then the universal covering 
map $f : \half \arrow U(r)$ is given by
\begin{displaymath}
	f(t)  = \exp(2\pi i \, \log(t) / \log(a) ).
\end{displaymath}
The deck group of $\half/U(r)$ is generated by $z \mapsto az$, and $A$ is the image of
the positive imaginary axis under $f$.  The values of $a$ and $r$ are
related by $r = f(-1) = \exp(-2\pi^2/\log(a))$.

\bold{Winding numbers.}
Consider now an immersed hyperbolic geodesic $\gamma \subset U(r)$ with both endpoints on the
unit circle $S^1$ (see Figure \ref{fig:annulus}).  
Any such geodesic in $U(r)$ can be presented as the image of an embedded geodesic
$\gammat \subset \half$ with endpoints $0<x_1<x_2$ on the positive real axis.
The {\em winding number} of $\gamma$ is the unique real number $N(\gamma) > 0$ such that
\begin{displaymath}
	x_2 = a^{N(\gamma)} x_1 .
\end{displaymath}
It is easily seen that $\gamma$ winds $N(\gamma)$ times around
the annulus $U(r)$, in the sense that the angle $\arg(z)$ increases strictly monotonically by 
$2\pi N(\gamma)$ as $z$ moves counter-clockwise from one end of $\gamma$ to the other.

The distance $d(z,A)$ is a strictly convex function on $\gamma$; in particular,
it assumes its minimum value at a unique point $p$.  The more $\gamma$ winds around the annulus,
the closer it comes to its core geodesic $A$; that is,
the minimum distance satisfies
\begin{equation}
\label{eq:Nd}
	N(\gamma) > N(\gamma') \implies d(\gamma,A) < d(\gamma',A) .
\end{equation}
This follows immediately from the fact that the hyperbolic distance from
$\gammat$ to the imaginary axis is a decreasing function of $x_2/x_1$.

\bold{Self--intersections.}
When the winding number $N(\gamma)$ exceeds $1$, the geodesic $\gamma$
has 1 or more self--intersections.
(To see this, just observe that the endpoints
$(ax_1, ax_2)$ of $a\gammat$ are linked with the endpoints
$(x_1,x_2)$ of $\gammat$ when $x_2 > a$.)
In particular, $\gamma$ has a unique
self--intersection point $q$ which is closest to $A$.
The point $q$ is opposite to $p$
and cuts off an embedded loop $\gamma^0 \subset \gamma$ encircling $A$.
The region between $A$ and $\gamma^0$ is a convex annulus in the
hyperbolic metric on $U(r)$.

If we remove $\gamma^0$ from $\gamma$, we are left with a path
$\gamma^- \subset U(r)$ whose geodesic representative has winding number $N(\gamma)-1$.

\makefig{The tower of coverings $\half \arrow \Mt(\lambda) \arrow M(\lambda)$.}{fig:covers}{
	\includegraphics[height=2.3in]{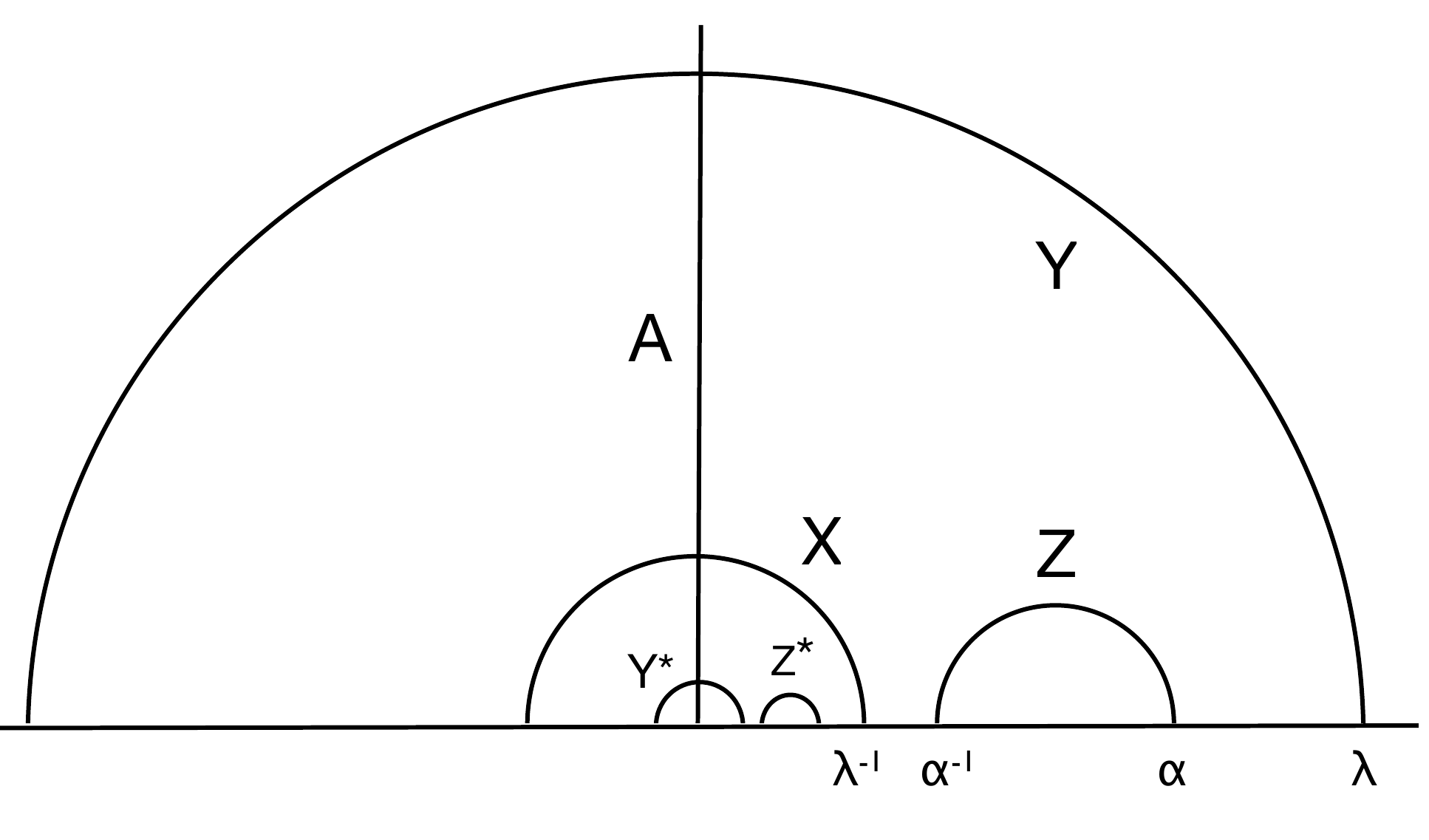}\\
	$\big \downarrow$\\
	\vspace{.1in}
	\includegraphics[height=1.9in]{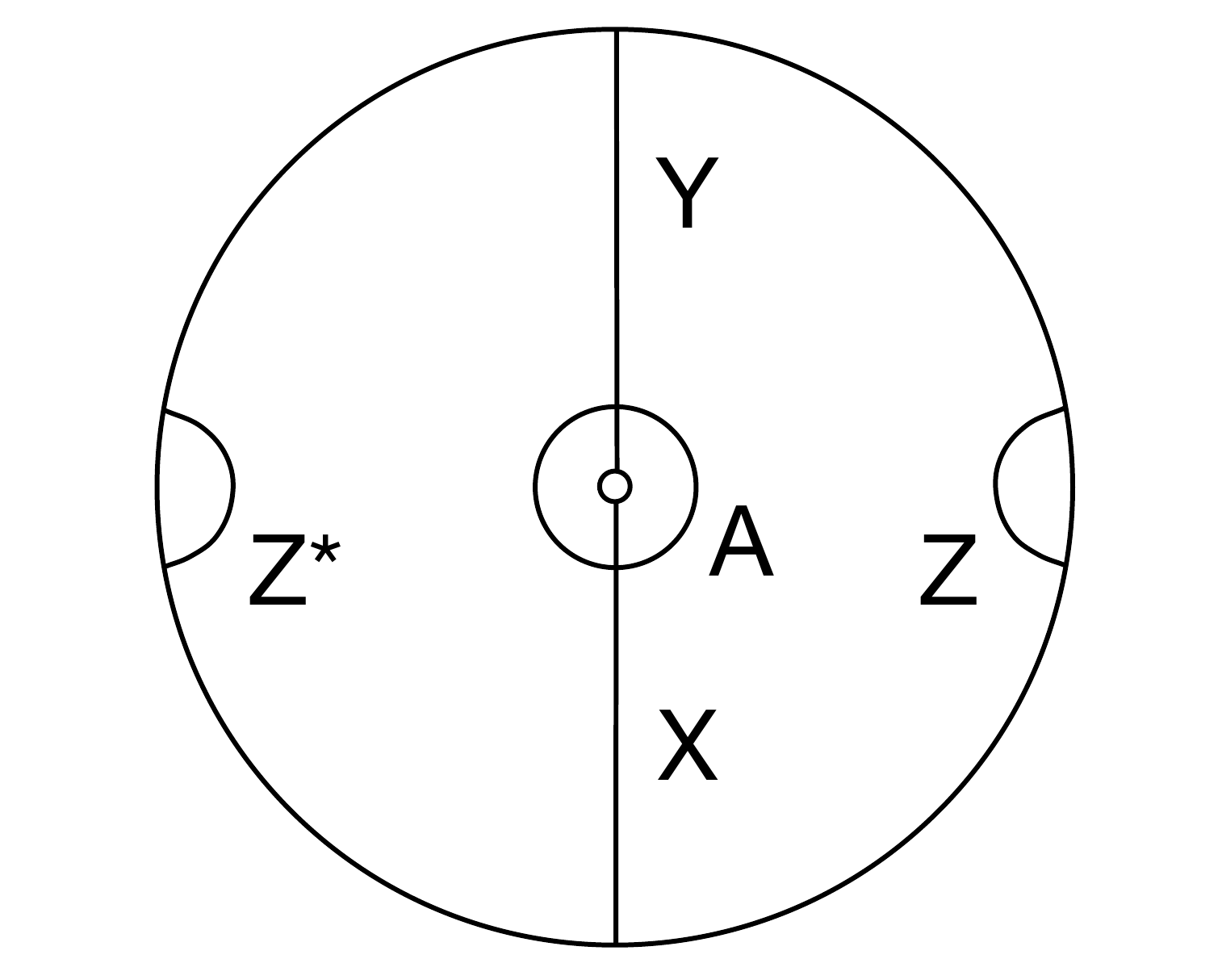}\\
	\vspace{.1in}
	$\big \downarrow$\\
	\vspace{.1in}
	\includegraphics[height=1.9in]{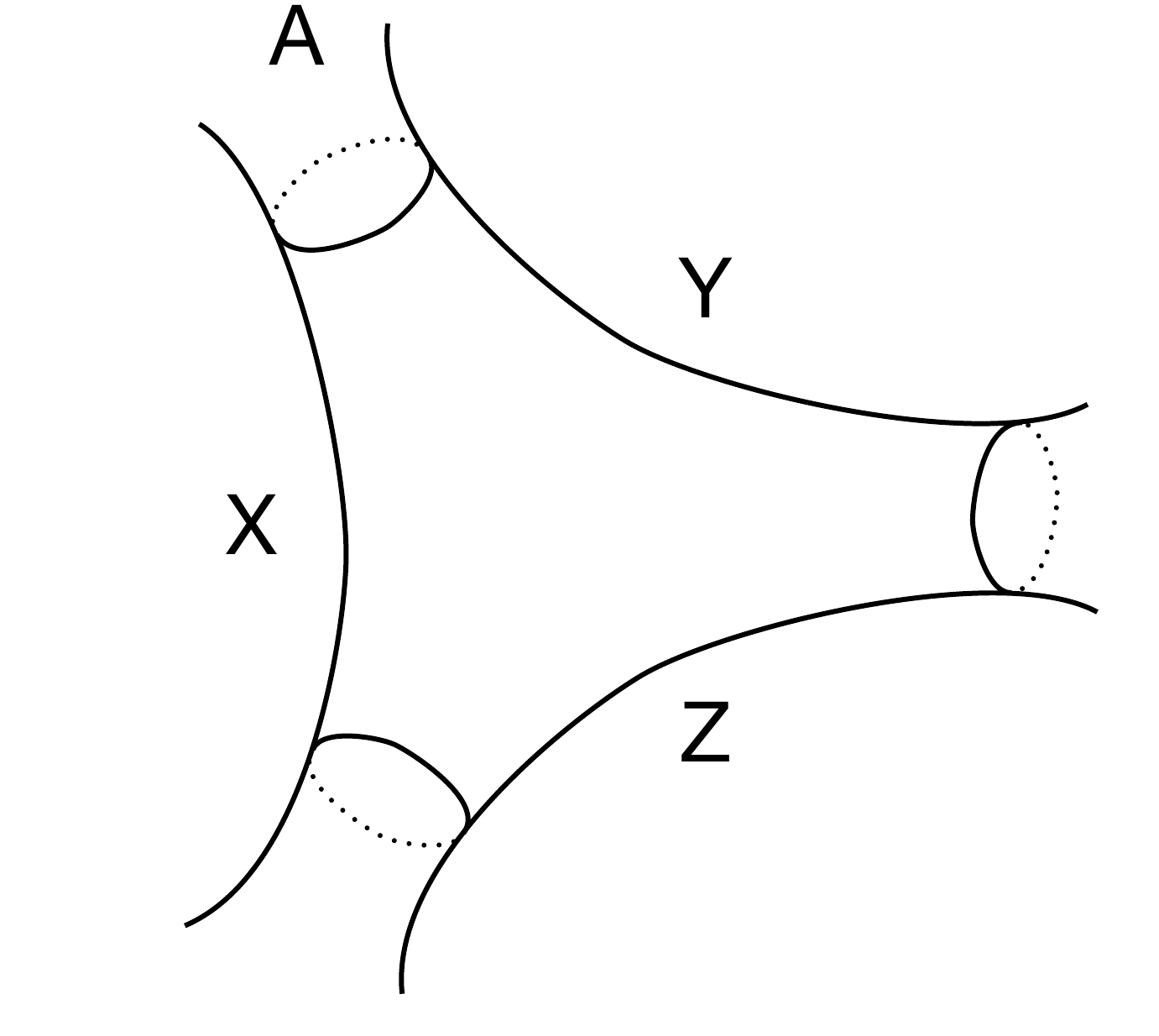}
}

\bold{Symmetric pairs of pants.}
To study the triply--punctured sphere $M = \half/\Gamma(2)$, it turns out to be useful to
consider instead the symmetric pair of pants $M(\lambda)$ with cuffs of length
$4 \log \lambda$.  
(It is well known that there is a unique pair of pants with cuffs of given
lengths; see e.g. \cite[\S 3.1.5]{Imayoshi:Taniguchi:book}.)

Given $\lambda > 1$, a concrete model for the surface $M(\lambda)$ can be
constructed as follows.
First, consider the hyperbolic geodesics $X, Y \subset \half$ defined by
$|t| = \lambda^{-1}$ and $|t|=\lambda$ respectively.  Let $Z \subset \half$ be the unique
geodesic disjoint from these two, and resting on the positive real axis, such that
\begin{displaymath}
	d(X,Y) = d(Y,Z) = d(Z,X) .
\end{displaymath}
The endpoints of $Z$ are given by $\alpha$ and $\alpha^{-1}$, where
$1 < \alpha < \lambda$.
The value of $\alpha$ is uniquely determined by $\lambda$ and satisfies
\begin{equation}
\label{eq:alpha}
	\alpha(\lambda) \asyto 1 + 1/\lambda
\end{equation}
as $\lambda \arrow \infty$.  It can be computed using cross-ratios and
the fact that there is a M\"obius transformation that cyclically
permutes $X$, $Y$ and $Z$.

Letting $x,y,z \mem \Isom(\half)$ denote the reflections through the geodesics $X, Y$ and $Z$ respectively,
we obtain a natural isometric action of the group
\begin{displaymath}
	\Gt = \{x,y,z \st x^2 = y^2 = z^2 = \id \} 
\end{displaymath}
introduced in \S\ref{sec:background}.
The orientation--preserving subgroup $G \subset \Gt$ is generated by $xy$, $yz$ and $zx$,
and the quotient space
\begin{displaymath}
	M(\lambda) = \half/G
\end{displaymath}
is the desired symmetric pair of pants.
A fundamental domain for the action of $G$ on $\half$ is given by
the region bounded by $X,X^*$ and $Z,Z^*$, where $X^*$ and $Z^*$ are the images of $X$ and $Z$ under reflection through $Y$
(see the top of Figure \ref{fig:covers}).

\bold{Intersection numbers.}
Note that as $\lambda \arrow 1$, the infinite volume surface $M(\lambda)$ converges geometrically to the finite volume
surface $M = \half/\Gamma(2)$, which we can regard as $M(1)$.

The self--intersection number $I(w)$ of a word $w \mem G$ can be
computed using its geodesic representative on $M(\lambda)$ just
as well as on $M(1)$; in fact, geodesics for any hyperbolic metric
automatically minimize intersection numbers (see e.g. \cite[Expos\'e 3, Theorem 15]{FLP}).
It turns out that the most advantageous geometry for our purposes arises,
not when $\lambda=1$, but when $\lambda$ is large.

\bold{The annular covering space.}
Note that $xy \mem G$ acts on $\half$ by the hyperbolic transformation
$t \mapsto at$, where $a  = \lambda^4$.  This transformation stabilizes the
vertical geodesic $A = i \reals_+ \subset \half$.

To analyze runs of type $xy$, we will use the intermediate covering space
\begin{displaymath}
	\Mt(\lambda) = \half / \brackets{xy} \isom U(r),
\end{displaymath}
where $r = \exp(-2\pi^2/\log a)$ as before.  The tower of coverings
\begin{displaymath}
	\half \arrow
	\Mt(\lambda) \arrow 
	M(\lambda) 
\end{displaymath}
is shown in Figure \ref{fig:covers}.  
The geodesics $A,X,Y,Z$ and $Z^*$ in $\half$ map to geodesics in $\Mt(\lambda)$ which for simplicity we
will continue to denote by the same letters.
Note that $A$ descends to the core geodesic of $\Mt(\lambda)$,
while $X, Y, Z$ and $Z^*$ map to embedded geodesics cutting $\Mt(\lambda)$ into disks.

We now fix any value of $\lambda$ such that
\begin{displaymath}
	\alpha = \alpha(\lambda) < \lambda^{1/2} .
\end{displaymath}
(Such a value exists by (\ref{eq:alpha}); in fact, any $\lambda \ge 2$ will do.)
This choice insures that $Z$ and $Z^*$ each cut off an arc of $S^1$ of length less than $\pi/2$.
As a useful consequence, we can estimate the winding number of a geodesic $\gamma \subset \Mt(\lambda)$
by combinatorial means.

\begin{lemma}
\label{lem:N}
Suppose a geodesic $\gamma \subset \Mt(\lambda) \isom U(r)$ first crosses $Z \cup Z^*$, then crosses $X \cup Y$ a total of $\ell$ times,
and finally crosses $Z \cup Z^*$ again.  Then its winding number satisfies
\begin{displaymath}
	\left| N(\gamma) - \frac{\ell}{2} \right| < \frac{1}{4} \cdot
\end{displaymath}
\end{lemma}

\bold{Proof.}
Suppose for simplicity that $\gamma$ first crosses $Z$ and winds counterclockwise around $U(r)$.
Then $\gamma$ has a lift to a geodesic $\gammat \subset \half$ 
that begins at a point $x_1$ in the interval bounded by the endpoints $(\alpha^{-1},\alpha)$ of $Z$,
and ends at a point $x_2$ in the interval bounded by the endpoints of $\lambda^{2\ell} Z$.
It follows that $\lambda^{2\ell}/\alpha^2 \le x_2/x_1 \le \lambda^{2\ell} \alpha^2$.   
Since $N(\gamma) = \log(x_2/x_1)/(4\log\lambda)$ and $1 < \alpha^2 < \lambda$,
the result follows.
\qed

An example of $\gamma \subset \Mt(\lambda)$ with $\ell = 3$ is shown in Figure \ref{fig:wind3}.

\makefig{The winding number of a geodesic $\gamma$ is estimated by half 
the number of times it crosses $X \cup Y$.}{fig:wind3}{
	\includegraphics[height=2.3in]{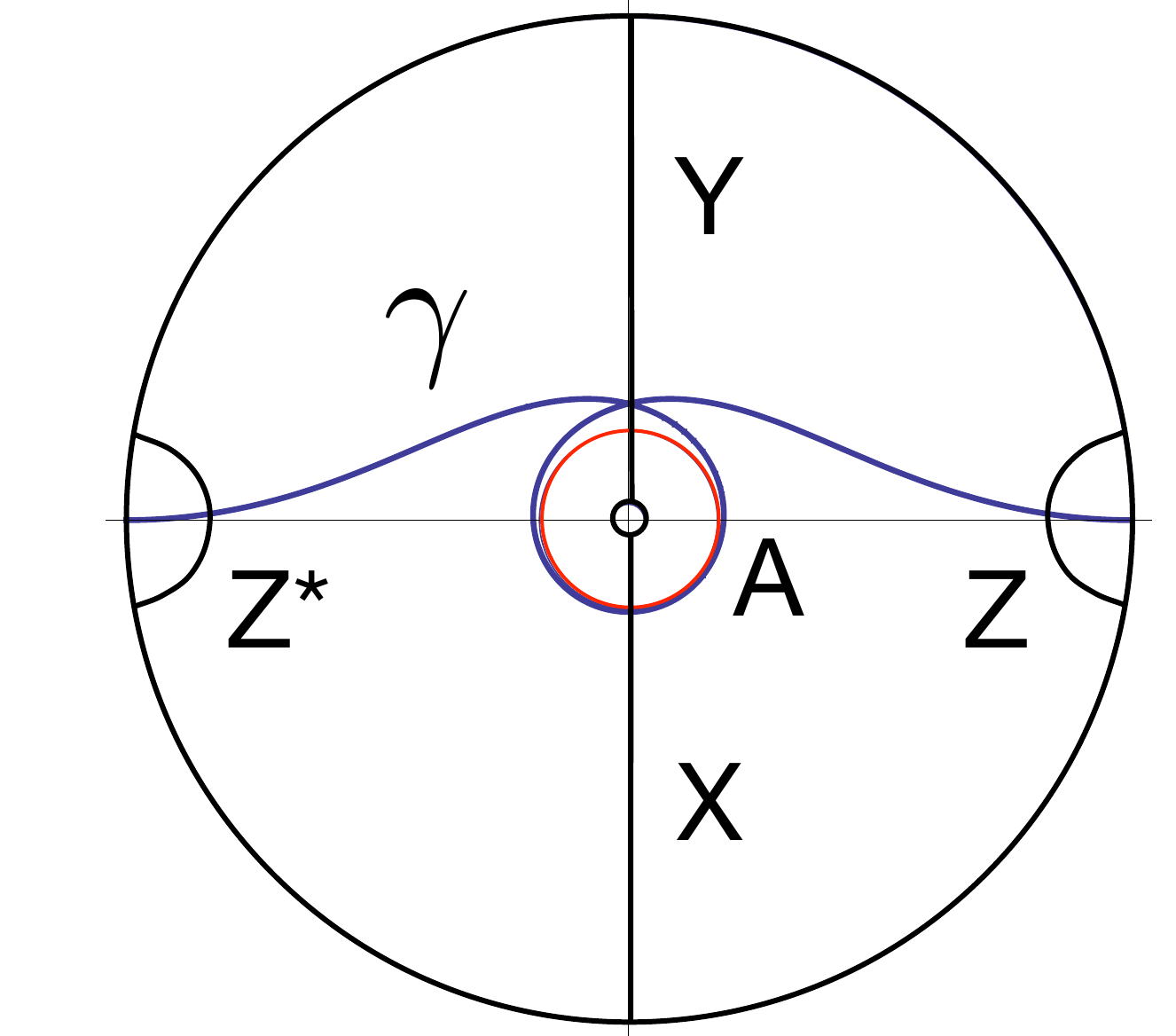}
}

\bold{Proof of Theorems \ref{thm:int1}, \ref{thm:int2} and \ref{thm:int3}.}
For concreteness we will analyze runs of type $xy$; the same argument applies to the other two types of runs with 
only notation modifications.

Given a class $[w] \mem \cG$, we can choose a representative cyclically reduced word in this class whose first letter is $z$.
Then we have
\begin{displaymath}
	w = z r_1 z r_2 \ldots z r_n,
\end{displaymath}
where the $r_i$ are maximal subwords in which no $z$ appears.  The $r_i$ of length two or more
are exactly the $xy$ runs of $w$.

There is a corresponding factorization $\gamma  = \gamma_1 * \cdots *\gamma_n$
of the closed geodesic $\gamma \subset M(\lambda)$ representing $[w]$, where the endpoints of each segment $\gamma_i$ lie on
$Z$ and the interior meets $X \cup Y$ a total of $\ell(r_i)$ times.

Let $F \subset \Mt(\lambda)$ be the region bounded by $Z \cup Z^*$ (see Figure \ref{fig:covers}).  The surface $M(\lambda)$
is obtained from $F$ by gluing its two geodesic boundary components together.  Thus we can regard each segment 
in the factorization of $\gamma$ as an arc $\gamma_i \subset F$ with endpoints of $Z \cup Z^*$.  

Let $\gammat_i \subset \Mt(\lambda)$ denote the extension of $\gamma_i$ to a complete geodesic in the annulus.
By Lemma \ref{lem:N}, the winding number of this geodesic is well--approximated by $\ell(r_i)/2$; that is, we have
\begin{equation}
\label{eq:N}
	\left| N(\gammat_i) - \frac{\ell(r_i)}{2} \right| < 1/4 .
\end{equation}

Let $I$ denote the set of indices such that $\ell(r_i) \ge 3$; these are the $xy$--runs of $w$ where contraction is possible.
For each $i \mem I$, we have $N(\gammat_i) > 1$, and hence $\gamma_i$ has a self--intersection that cuts off an
innermost loop $\gamma_i^0$ encircling the core geodesic $A \subset \Mt(\lambda)$.
Let $\gamma_i^-$ denote $\gamma_i$ with the loop $\gamma_i^0$ removed.  Then the loop
\begin{displaymath}
	\gamma^- = \gamma_1 * \cdots * \gamma_i^- * \cdots * \gamma_n
\end{displaymath}
represents the homotopy class $[w^-[r_i]]$.  We now distinguish two cases.

(1) Suppose $\gamma_i^0 \cap \gamma_j \neq \emptyset$ for some $j \neq i$.  Then $\gamma$ has three more self--intersections
that $\gamma^-$.  Indeed, the arc $\gamma_j$ must enter and exit the annulus bounded by $A \cup \gamma_i^0$, resulting in
two self--intersections; and the loop $\gamma_i^0$ comes from an intersection of $\gamma_i$ with itself.
Thus:
\begin{displaymath}
	I(w^-[r]) \le I(w) - 3.
\end{displaymath}
In particular, $r_i$ is not exceptional.

(2) Now suppose $\gamma_i^0 \cap \gamma_j = \emptyset$ for all $j \neq i$.  Then $\gamma$ still has one more self--intersection
than $\gamma^-$, and hence
\begin{equation}
\label{eq:1}
	I(w^-[r]) \le I(w) - 1.
\end{equation}
Moreover, since $\gamma$ was a geodesic, $\gamma^-$ admits a regular homotopy to its geodesic
representative; thus $I(w^-[r]) = I(w)-1 \mod 2$.   So in this case, either $I(w^-[r])=I(w)-1$ --- and $r$ is exceptional;
or we have $I(w^-[r]) \le I(w) -3$, just like in case (1).

For case (2) to hold, we must have
\begin{equation}
\label{eq:close}
	d(\gamma_i^0,A) < d(\gamma_j,A)
\end{equation}
for all $j \neq i$.  This shows that (2) holds for at most one value of $i$, and hence $w$ has at most  
one exceptional run of type $xy$ (Theorem \ref{thm:int1}).  Moreover, (\ref{eq:close}) implies, by
relation (\ref{eq:Nd}), that $N(\gammat_i) > N(\gammat_j)$ for all $i \neq j$, and hence
$\ell(r_i) \ge \ell(r_j)$ for all $i \neq j$ by equation (\ref{eq:N}).  This completes the proof of Theorem \ref{thm:int2}.

Finally we prove Theorem \ref{thm:int3}.  Suppose $\ell(r_i) > \ell(r_j)$ for all $j \neq i$.  
Then, by equations (\ref{eq:N}) and (\ref{eq:Nd}), we have
\begin{displaymath}
	d(\gamma_i,A) < d(\gamma_j,A)
\end{displaymath}
for all $j \neq i$.  Let $p \mem \gamma_i$ be the point closest to $A$.  By attaching a new loop to $\gamma_i$ at $p$
that runs once around $A$, we introduce one new self--intersection and obtain a path
\begin{displaymath}
	\gamma^+ = \gamma_1 * \cdots * \gamma_i^+ * \cdots *\gamma_n
\end{displaymath}
representing the class $[w^+[r_i]] \mem \cG$.  This shows that
\begin{displaymath}
	I(w^+[r]) \le I(w) + 1 .
\end{displaymath}
On the other hand, if we erase this loop then the self--intersection number goes down by at least one, by Theorem \ref{thm:int2},
and hence equality holds.
\qed

\section{Length and self--intersections}
\label{sec:minusone}

With the surgery bounds in place, it is now straightforward to establish a lower bound on the defect
\begin{displaymath}
	\delta(w) = I(w) - L(w) .
\end{displaymath}

\begin{theorem}
\label{thm:short}
For any cyclically reduced word $w \mem G$, we have $\delta(w) \ge -1$.
\end{theorem}

\bold{Multiple loops.}
As a preliminary, we remark that any primitive, cyclically reduced
word $w \mem G$ satisfies
\begin{equation}
\label{eq:wn}
	I(w^n) = 
	n^2 I(w) + n - 1 
	\AND
	L(w^n) = nL(w).
\end{equation}
The formula for $L(w^n)$ is immediate, and the formula
for $I(w^n)$ is well known (see e.g. \cite[Theorem 6]{deGraaf:Schrijver:crossings}).
In fact, an optimal loop representing
$[w^n]$ can be obtained from $n$ parallel copies of an optimal
loop for $[w]$ by cyclically braiding the strands to form a single loop.
The parallel copies yield $n^2 I(w)$ crossings, and the braiding
accounts for $n-1$ more.

\bold{Proof of Theorem \ref{thm:short}.}
We first remark that if this bound holds for a primitive, cyclically reduced word $w \mem G$, then it also holds for $w^n$, $n>1$.
Indeed, if we know that $I(w) \ge L(w) - 1$, then by equation
(\ref{eq:wn}) we also have:
\begin{displaymath}
	I(w^n) = n^2 I(w) + n - 1 \ge n^2 L(w) - n^2 + n - 1 \ge n L(w) - 1 = L(w^n) - 1,
\end{displaymath}
and hence $\delta(w^n) \ge -1$ as well.

We can now argue by contradiction.  Suppose the bound is false.
Let $w \mem G$ be a word of minimal length with $\delta(w) < -1$.
By the preceding remark, $w$ is primitive.
Suppose $w$ has a run $r$ of length three or more, and let $w' = w^-[r]$.  Then $L(w') = L(w)-1$ and, by Theorem \ref{thm:int2}, we have
$I(w') \le I(w) -1$; hence $\delta(w') \le \delta(w) < -1$, contradicting the fact that $w$ has minimal length.

Thus $w$ has no run of length 3.  But there are only two types of words in $G$ with this property:
one is $w=xy$, and the other is $w=xyzxyz$.   The first has $\delta(w) = -1$, and the second has $\delta(w) = 0$.
\qed

\section{Motifs}
\label{sec:motifs}

In this section we will define the set of motifs
$\cM \subset \cG$, and show:

\begin{theorem}
\label{thm:bin}
The number of closed geodesics of length $L$ and defect $\delta$ is given by:
\begin{displaymath}
	N_\delta(L) = \sum_{[w] \mem \cM \st \delta(w) = \delta} \binom{L - L(w) + \rho(w) - 1}{\rho(w) - 1} \cdot
\end{displaymath}
\end{theorem}

\bold{Rank and descendants.}
Let $w$ be a cyclically reduced word with $[w] \mem \cG$.
Let us say a run $r$ of $w$ is {\em distinguished} if $\ell(r) > \ell(s)$ for every other run $s$ of the same type as $r$.

\begin{prop}
\label{prop:delta}
For any distinguished run $r$ of $w$, we have $\delta(w^+[r]) = \delta(w)$.
\end{prop}

\bold{Proof.}
By Theorem \ref{thm:int3}, expansion of a distinguished run increases both $L(w)$ and $I(w)$ by one.
\qed

The {\em rank} $\rho(w)$ is the number of its distinguished runs.
Since there is at most one distinguished run of each type, we have $0 \le \rho(w) \le 3$.

The set of {\em descendants} $\cD(w) \subset \cG$ is defined to be the smallest set containing 
$[w]$ and closed under expansion, in the sense that
\begin{displaymath}
	[v] \mem \cD(w) \implies [v^+[r]] \mem \cD(w)\;\text{for every distinguished run $r$ of $v$.}  
\end{displaymath}
For example, if $w$ has rank three, then it has three runs that can be independently extended, and thus
its set of descendants has the form
\begin{displaymath}
	\cD(w) = \{ [w_1 (xy)^i w_2 (yz)^j w_3 (zx)^k] \st i,j,k \ge 0 \} .
\end{displaymath}
In this case the number of descendants of $[w]$ of length $L$ is the same as the number of ways to
express $L$ in the form $L(w) + i + j + k$ with $i,j,k \ge 0$.  More generally, keeping in mind our conventions 
(\ref{eq:bin1}) and (\ref{eq:bin2}) on  binomial coefficients, we have:

\begin{prop}
\label{prop:binom}
The number of descendants of $w$ of length $L$ is given by 
\begin{displaymath}
	|\{[v] \mem \cD(w) \st L(v) = L\}| = 
	\binom{L-L(w) + \rho(w) - 1}{\rho(w) - 1} .
\end{displaymath}
\end{prop}

\bold{Motifs.}
A word $w$ with $[w] \mem \cG$ is a {\em motif} if, for any run 
$r$ of $w$ with $\ell(r) \ge 4$, there is another run $s$ of the same
type with $\ell(r) \le \ell(s)+2$.  
In a motif, there is therefore either a tie or a near--tie for the longest
run of a given type, unless there is a unique run $r$ of that
type and $\ell(r) \le 3$.

Let $\cM \subset \cG$ denote the set of all $[w]$ such that $w$ is a motif.

\begin{prop}
\label{prop:motif}
Every $[v] \mem \cG$ is a descendant of a unique motif $[w]$.  
\end{prop}

\bold{Proof.}
Given any $[v] \mem \cG$, repeatedly contract the distinguished runs $r$ of $v$ until either $\ell(r) \le 3$ or
$\ell(r) \le \ell(s)+2$ for some run $s$ of the same type as $r$.  The result is a motif $[w] = f([v])$
with $[v] \mem \cD(w)$.    Since we only contract when $\ell(r) \ge 4$, runs other than $r$ are unaffected
at each step, and thus the order of contraction does not change the result.
It is now readily verified by induction on $L(v)$ that $f([v]) = [w]$ for all motifs $[w]$ and
all $[v] \mem \cD(w)$; hence different motifs have different descendants.
\qed

\bold{Proof of Theorem \ref{thm:bin}.}
By Propositions \ref{prop:delta} and \ref{prop:motif}, the set of $[v] \mem \cG$ with $\delta(v) = \delta$
is the disjoint union of the descendants of the motifs $[w]$ with $\delta(w) = \delta$; now apply \ref{prop:binom}.
\qed

\bold{Examples.}
The 27 motifs $[w] \mem \cG$ with $\delta(w)= 0$ are listed in Table \ref{tab:motifs}.
For brevity we have chosen only one representative motif from each orbit of $\Aut(G)$
acting on $\cM$.  (Automorphisms of $G$ include, for example, permutations of the generators $x,y,z$.)  
The size of the orbit is denoted by $C(w)$, and a representative from the orbit is given in the final column.
These 27 motifs correspond to the loops that appear in Figure \ref{fig:N0} for $i,j,k \mem \{-1,0,1\}$.

\maketab{
	The 27 motifs with $\delta(w)=0$ come from 6 patterns.
	}{tab:motifs}{
\begin{tabular}{|c|c|c|c|}
	\hline
	$\rho(w)$ & $L(w)$ & $C(w)$ & $w$ \\ \hline
	\hline
  0 & 3 & 1 & $x . y . z . x . y . z$\\ \hline
  1 & 4 & 6 & $x . y . x . y . z . x . y . z$\\ \hline
  2 & 5 & 6 & $x . y . x . y . z . x . y . z . x . z$\\ \hline
  2 & 5 & 6 & $x . y . x . y . z . x . y . z . y . z$\\ \hline
  3 & 6 & 6 & $x . y . x . y . z . x . y . z . y . z . x . z$\\ \hline
  3 & 6 & 2 & $x . y . x . y . z . x . z . x . y . z . y . z$\\ \hline
\end{tabular}
}

The 27 motifs with $\delta=0$ yield formula (\ref{eq:N0}) for $N_0(L)$.
Similarly, the 12 motifs with $\delta = -1$ yield the formula:
\begin{displaymath}
	N_{-1}(L) = 3 \binom{L-2}{1} + 3 \binom{L-1}{1} + 6 \binom{L-2}{2},
\end{displaymath}
and the 153 motifs with $\delta=1$ yield the formula:
{\small
\begin{displaymath}
	N_1(L) = 
	 3 \binom{L-5}{-1}+
	 24 \binom{L-5}{0}+
	 54 \binom{L-5}{1}+
	 12 \binom{L-4}{1}+
	 36 \binom{L-5}{2}+
	 24 \binom{L-4}{2}  \cdot
\end{displaymath}
}
There are 135 motifs with $\delta=2$,
603 with $\delta=3$, 564 with $\delta=4$, and $2391$ with $\delta=5$.

\section{Lengths of motifs}
\label{sec:rank}

In this section we will show that the length of a motif is controlled by its self--intersection number.
More precisely, we will prove: 

\begin{theorem}
\label{thm:rank}
For any motif $w$, we have
\begin{displaymath}
	\delta(w)+\rho(w)+3 \ge L(w) .
\end{displaymath}
\end{theorem}

We then complete our main objective, the proof of
Theorem \ref{thm:main}.
We will use the bound above to show that $N_\delta(L)$ is a
{\em finite} sum of binomial coefficients, and that
the equality $N_\delta(L) = p_\delta(L)$ holds for all $L \ge \delta+4$.

\bold{Thin motifs.}
A motif $w$ is {\em thin} if it has no run with $\ell(r) \ge 4$.
The idea of the proof of Theorem \ref{thm:rank} is to use contractions
to reduce to the case of thin motifs.  The following quadratic bound 
plays an important role:

\begin{theorem}
\label{thm:thin}
For any thin motif, we have $I(w) \ge L(w)^2/6$.
\end{theorem}

\bold{Proof.}
Under the change of variables $(a=xy,b=yz,c=zx)$, a thin motif becomes
a word in $(a,b,c)$ with no repeated letters, so the
result follows from Theorem \ref{thm:quad} in the Appendix.
\qed

\bold{Properties of thin motifs.}
The conclusion of Theorem \ref{thm:rank} is equivalent to the bound:
\begin{equation}
\label{eq:Delta}
	\Delta(w) + \rho(w) \ge -3 ,
\end{equation}
where
\begin{displaymath}
	\Delta(w) = I(w) - 2L(w) .
\end{displaymath}
The proof is essentially by induction on the length of $w$.
To verify the inequality for motifs of small length, we use the algorithm
described in the Appendix to explicitly compute $I(w)$ in many cases.
Here are the base cases to be used in the induction.

\begin{theorem}
\label{thm:basis}
The inequality $\Delta(w) + \rho(w) \ge -3$ holds for:
\begin{enumerate}
	\item[(a)]
All thin motifs; and
	\item[(b)]
All motifs with $L(w) \le 8$.
\end{enumerate}
The sharper inequality $\Delta(w) \ge -3$ holds for:
\begin{enumerate}
	\item[(c)]
All thin motifs with $L(w) \ge 6$, and
	\item[(d)]
All imprimitive words that involve all three generators $x,y,z \mem \Gt$.
\end{enumerate}
\end{theorem}

\bold{Proof.}
Statement (b) is verified by directly calculating $\Delta(w)$ for the 2904 motifs with $L(w) \le 8$,
using the formula for $I(w)$ given in Theorem \ref{thm:I} below.
A similar calculation gives $\Delta(w) \ge -3$ for the 536 thin motifs with $6 \le L(w) \le 9$.
For $L(w) \ge 10$, Theorem \ref{thm:thin} implies that any thin motif with $L(w) \ge 10$ satisfies $\Delta(w) \ge 10^2/6 - 20 > -4$;
hence $\Delta(w) \ge -3$, and (c) follows.  Clearly (c) and (b) imply (a).

To prove (d), note that if $w$ is a {\em primitive} word involving all three generators $x,y,z$, then $L(w) \ge 2$, 
and $I(w) \ge L(w)-1$ by Theorem \ref{thm:minusone}; while for $n \ge 2$, 
we have $L(w^n) = nL(w)$ and $I(w^n) = n^2 I(w) + n - 1$
by equation (\ref{eq:wn}), and therefore
\begin{eqnarray*}
	\Delta(w^n) & \ge & (n^2 - 2n) L(w) - n^2 + n - 1 \\
	& \ge & n^2 - 3n - 1 \ge -3 .
\end{eqnarray*}
\qed

\bold{Remarks.}
The thin motif $w = (xyxz)^2yz$ with $L(w)=5$ satisfies $\Delta(w)=-4$, so statement (c) cannot be sharpened.
The imprimitive word $w=(xy)^3$ also satisfies $\Delta(w) = -4$, so the requirement that $w$ involves all
three generators is necessary in statement (d).

One can avoid the roughly 3500 special cases implicit
in the proof of Theorem \ref{thm:basis} at the cost of
replacing the lower bound of $-3$ with, say, $-10$.
This would still suffice to prove Theorem \ref{thm:main}, but with
$\delta+4$ replaced by $\delta+11$.

\bold{Double contraction.}
To reduce Theorem \ref{thm:rank} to the case of thin motifs,
we will use the following observations about contractions, based on the results of \S\ref{sec:surgery}.

We will only contract runs with $\ell(r) \ge 4$.  This insures that the runs of $w^-[r]$, other than $r$ itself, remain the same as the
runs of $w$.  In particular, given two different runs $r,s$ of the same type, we can contract them in either order to obtain
the {\em double contraction} $w^-[r,s]$.
By Theorem \ref{thm:int2}, we have:
\begin{displaymath}
	\Delta(w^-[r]) \le \begin{cases}
			\Delta(r) + 1 & \text{if $r$ is exceptional},\\
			\Delta(r) - 1 & \text{if $r$ is not exceptional.}
			\end{cases}
\end{displaymath}
By Theorem \ref{thm:int1}, at most one of the runs $r$ and $s$ is exceptional, so we have:
\begin{displaymath}
	\Delta(w^-[r,s]) \le \Delta(w) .
\end{displaymath}
This observation is the crux of the argument.

\bold{Proof of Theorem \ref{thm:rank}.}
Suppose the desired bound is false, and let $w$ be a motif of minimum length with 
\begin{equation}
	\label{eq:viol}
	\Delta(w) + \rho(w) < - 3.
\end{equation}
Then $L(w) \ge 9$ by Theorem \ref{thm:basis} part (b), and any contraction $w'$ of $w$ with $\Delta(w') \le \Delta(w)$ must be primitive by part (d).

Consider the runs of $w$ of a given type, ordered so that $\ell(r_1) \ge \ell(r_2) \ge \cdots \ge \ell(r_n)$.  We can assume that
the exceptional run, if any, is $r_1$.  By the definition of a motif, either $\ell(r_1) \le 3$ or $\ell(r_1) \le \ell(r_2)+2$.

We claim that $\ell(r_i) \le 3$ for all $i \ge 3$; otherwise,
$w^-[r_i]$ would give a shorter word satisfying equation (\ref{eq:viol}).  
We also have $\ell(r_1) \le 5$; otherwise,
$w^-[r_1,r_2]$ would give a shorter word satisfying (\ref{eq:viol}).

We also have $\ell(r_2) \le 3$, for similar reasons.
That is, we cannot have $(\ell(r_1),\ell(r_2)) = (5,5), (4,4)$ or $(5,4)$.
Indeed, in the $(5,5)$ and $(4,4)$ cases, $w^-[r_2]$ would be 
a smaller motif satisfying (\ref{eq:viol}) 
(contracting $r_2$ increases $\rho(w)$ by one, but it also
decreases $\Delta(w)$ by one);
while in the $(5,4)$ case, $w^-[r_1,r_2]$ would be a smaller solution.

It follows that, for runs of a fixed type, either $\ell(r_1) \le 3$ or 
\begin{displaymath}
	\ell(r_2) < \ell(r_1) \le 5.
\end{displaymath}
In particular, the number of types with $\ell(r_1)>3$ is at most $\rho(w)$.  By contracting $r_1$ for each such type,
we obtain a thin motif $w'$.  Each contraction reduces the length of $w$ by one, and
increases $\Delta(w)$ by at most one; hence
\begin{displaymath}
	\Delta(w) + \rho(w) \ge \Delta(w'),
\end{displaymath}
and $L(w') \ge L(w) - \rho(w) \ge 9-3 = 6$.  Therefore $\Delta(w') \ge -3$ by
Theorem \ref{thm:basis}(c),  contradicting our assumption (\ref{eq:viol}).
\qed

\bold{Remark.}
The constant 3 in the statement of Theorem 
\ref{thm:rank} cannot be improved; for example, the 
motif $w = xyzxyz$ satisfies $\delta(w)=\rho(w)=0$ and $L(w) = 3$.

\bold{Proof of Theorem \ref{thm:main}.}
Fix $\delta \ge -1$.  By Theorem \ref{thm:rank}, any motif $[w] \mem \cM$ with
$\delta(w) = \delta$ satisfies
\begin{displaymath}
	L(w) \le \delta + \rho(w) + 3\le \delta + 6,
\end{displaymath}
since $\rho(w) \le 3$.
Thus there are only finitely many such motifs, and hence Theorem 
\ref{thm:binom} expresses $N_\delta(L)$ as a finite sum of binomial
coefficients.  Each binomial coefficient has the form
\begin{displaymath}
	b(L) = \binom{L - L(w) + \rho(w) - 1}{\rho(w) - 1}
\end{displaymath}
for some $[w] \mem \cM$.  Now recall that $\binom{n}{k}$ 
agrees with a polynomial of degree $k$ in $n$ for all $n \ge 0$.
Theorem \ref{thm:rank} also implies that
\begin{displaymath}
	L - L(w) + \rho(w) - 1 \ge L - (\delta + 4).
\end{displaymath}
Thus for $L \ge \delta+4$, 
Theorem \ref{thm:binom} expresses the value of $N_\delta(L)$ as
a quadratic polynomial in $L$.
\qed

\appendix

\section{Appendix:  Counting self--intersections}
\label{sec:app}

\makefig{Generators $a,b,c$ for $\pi_1(M)$, and the corresponding Cayley graph $H$.}{fig:cayley}{
\includegraphics[height=1.9in]{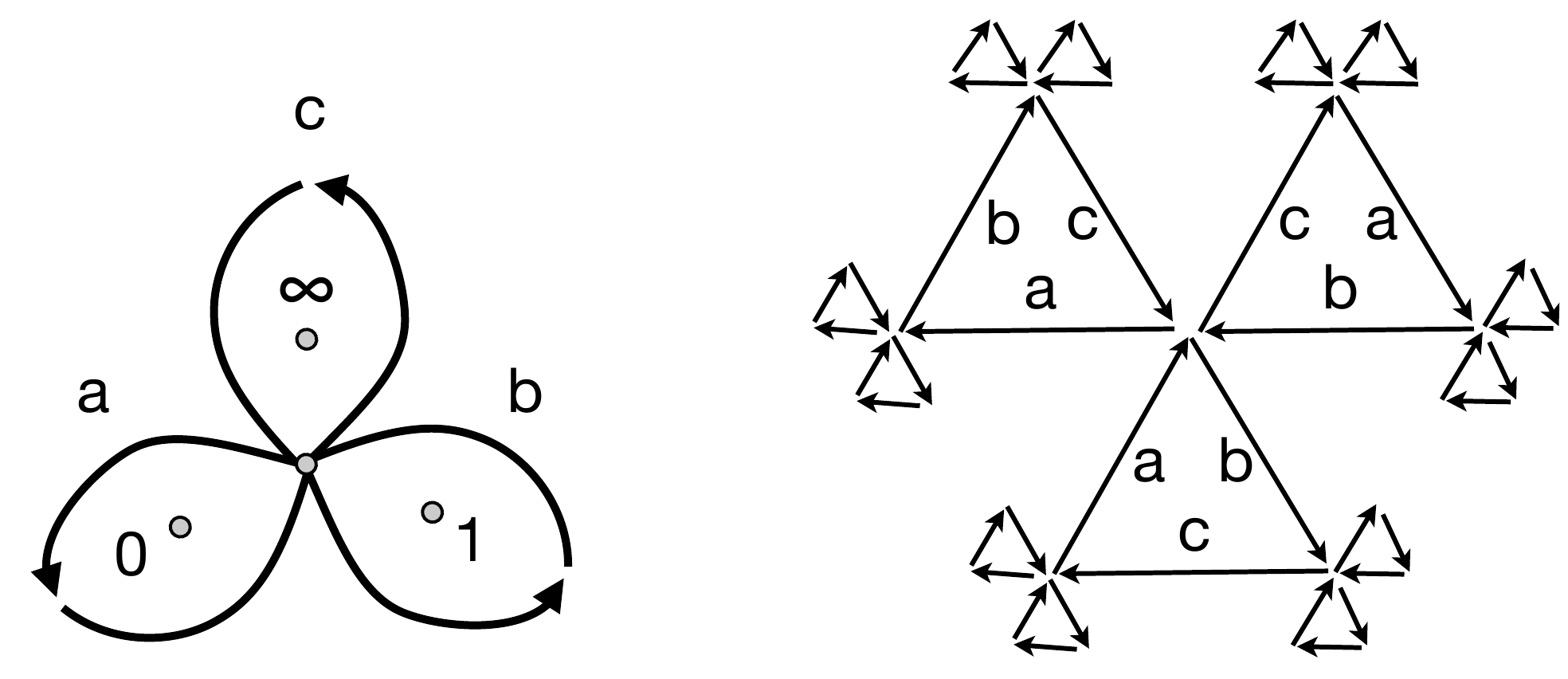}}

In this section we give a combinatorial formula for the self--intersection 
number $I(w)$.
This formula was used to produce the data in Table \ref{tab:main}.
We also obtain the following useful quadratic lower bound on the intersection number:

\begin{theorem}
\label{thm:quad}
Let $w \mem G$ be a cyclically reduced word in the generators $a,b,c$ 
and their inverses.  Suppose that $L(w) \ge 2$ and no two consecutive letters of $w$ are the same.  
Then we have $I(w) \ge L(w)^2/6$.
\end{theorem}

\bold{Group theory.}
Let $G = \pi_1(M,p)$, where $M = \chat - \{0,1,\infty\}$ and $p$ is a 
point in the upper halfplane.
In this section we will use the presentation
\begin{equation}
\label{eq:pres}
	G = \brackets{a,b,c \st abc=\id} 
\end{equation}
for $\pi_1(M,p)$.   It is related to the description of $G$ in \S\ref{sec:background} by the change of variables
\begin{equation}
\label{eq:abc:xyz}
	a = xy, \;\;
	b = yz \AND
	c = zx .
\end{equation}

The presentation (\ref{eq:pres}) corresponds to the generating loops for $\pi_1(M,p)$ shown at the left in Figure \ref{fig:cayley}.
A portion of the Cayley graph $H$ of $G$ is show at the right.  There is a natural realization of
$H$ as a graph in the plane, consistent with the ribbon graph structure on $H/G$ inherited from $M$.
The graph $H$ is a proxy for the universal cover $\half$ of $M$.  It carries a canonical metric
where each edge has length one.

\bold{Intersection numbers.}
Any two points in $H$ lie on a unique complete geodesic $\delta \subset H$.
The intersection number $i(\alpha,\beta)$ of a pair of geodesics in $H$ is defined to be
$1$ if $\alpha$ and $\beta$ {\em cross} in the plane; equivalently, if their endpoints are
linked at infinity.  Otherwise $i(\alpha,\beta) = 0$.
Note that geodesics can meet without crossing; in particular, we have $i(\alpha,\alpha)=0$.

\bold{Reduced words.}
Every $w \mem G$ has a unique expression as a {\em reduced word} in the generators $a,b,c$ and their inverses
$\abar, \bbar, \cbar$.
A reduced word has minimal length among all products representing the same element of $G$.
Concretely, this means that consecutive letters like $ab$, $bc$ and $ca$ must be avoided,
since they can be replaced by $\cbar, \abar$ and $\bbar$ respectively.

We let $L(w)$ denote the length of the $w$ as a reduced word in $a,b,c$.  This notion of length is consistent with the definition
$L(w) = \ell(w)/2$ given in \S\ref{sec:background}, as can be seen using equation (\ref{eq:abc:xyz}).

\bold{Stabilizers of geodesics.}
We say $w \mem G$ is {\em cyclically reduced} if it minimizes $L(w)$ among all elements in its conjugacy class.  
Concretely, this means that patterns like $ab$ are also be avoided when we regard the first and
last letters of $w$ as consecutive.

Aside from the identity, every element $w \mem G$ stabilizes a unique geodesic $\delta_w \subset H$.
This geodesic passes through the identity element $e \mem H$ if and only if $w$ is cyclically reduced.
If $w$ is primitive, then it generates the stabilizer of $\delta_w$.
Similarly, $w$ stabilizes a unique geodesic $\gammat_w$ in the universal cover $\half$ of $M$ provided
$L(w) \ge 2$.

\bold{Self--intersection numbers.}
Now let us fix a primitive, cyclically reduced word
\begin{displaymath}
	w = g_1 g_2 \cdots g_n \mem G
\end{displaymath}
with $n = L(w) \ge 2$.  Its cyclically reduced conjugates are given by
\begin{displaymath}
	w_i = g_i g_{i+1} \cdots g_n g_1 \cdots g_{i-1} ,
\end{displaymath}
$i=1,2,\ldots,n$.  Let $\delta_i$ denote the geodesic through the origin $e \mem H$ stabilized by $w_i$,
and let
\begin{equation}
\label{eq:deg}
	\deg_e(\delta_i \cup \delta_j) = |\{g_i,\gbar_{i-1},g_j,\gbar_{j-1}\}|
\end{equation}
denote the number of edges of $\delta_i \cup \delta_j$ incident to the vertex $e \mem H$.

Algorithms to compute $I(w) = I(\gamma_w)$ in various situations 
are well--known; see e.g.  \cite{Cohen:Lustig:sccs}, \cite{Despre:Lazarus:sccs}.
For the case at hand, we will show:

\begin{theorem}
\label{thm:I}
For any primitive, cyclically reduced word $w \mem G$, we have
\begin{equation}
\label{eq:I1}
	I(w) = \frac{1}{4} \sum_{1 \le i,j \le n} i(\delta_i,\delta_j) (\deg_e(\delta_i,\delta_j) - 2) .
\end{equation}
\end{theorem}

\bold{Proof.}
We take as our point of departure the following expression for 
the self--intersection number as a sum over double cosets:
\begin{equation}
\label{eq:I2}
	I(w) = \frac{1}{2}
		\sum_{[g] \mem \brackets{w} \bs G / \brackets{w}} i(g \cdot \delta_w,\delta_w) .
\end{equation}
This formula is easily justified using hyperbolic geometry, by observing that
$g \cdot \gammat_w$ intersects $\gammat_w$ in $\half$ iff $i(g \cdot \delta_w,\delta_w) = 1$.  Each nonzero term in the sum corresponds to a multiple
point $p \mem \gamma$ together with an ordered pair of branches
of $\gamma$ through $p$, and contributes $1/2$ to total intersection
number $I(\gamma_w) = I(w)$.

To connect formulas (\ref{eq:I2}) and (\ref{eq:I1}), suppose that
$i(g \cdot \delta_w,\delta_w) = 1$.  Then 
\begin{displaymath}
	(g\cdot \delta_w) \cap \delta_w 
\end{displaymath}
is a compact geodesic interval in $H$ with endpoints $h_1, h_2 \mem G$.  Translating
by $h_1^{-1}$, we obtain a pair of geodesics 
passing through $e$, each stabilized by a conjugate of $w$; hence
\begin{displaymath}
	(h_1^{-1} g \cdot \delta_w, h_1^{-1} \cdot \delta_w) = (\delta_i,\delta_j)
\end{displaymath}
for some $1\le i,j \le n$.  The intersection number is preserved by
translation, so $i(\delta_i,\delta_j) = 1$.

If $h_1 = h_2$ then $\delta_i$ and $\delta_j$ cross transversally
at $e$; hence $\deg_e(\delta_i \cap \delta_j) = 4$ and we
obtain a contribution of $1/2$ to the sum in formula (\ref{eq:I1}).
If $h_1 \neq h_2$ then the degree is 3, and we obtain two terms,
one from $h_1$ and one from $h_2$,
each contributing $1/4$ to the same sum.  
All the nonzero terms in equation (\ref{eq:I1}) arise in this
way, and hence the two sums are equal.
\qed

\bold{Algorithmic considerations.}
The term $i(\delta_i,\delta_j)$ appearing in 
formula (\ref{eq:I1}) can be readily computed
by comparing the cyclic orderings of the edges of $\delta_i$ and $\delta_j$
at the endpoints of the segment $\delta_i \cap \delta_j$.
The degree is computed by equation (\ref{eq:deg}).

\makefig{Chords represent geodesics passing through the origin $e \mem H$.}{fig:abc}{
\includegraphics[height=1.0in]{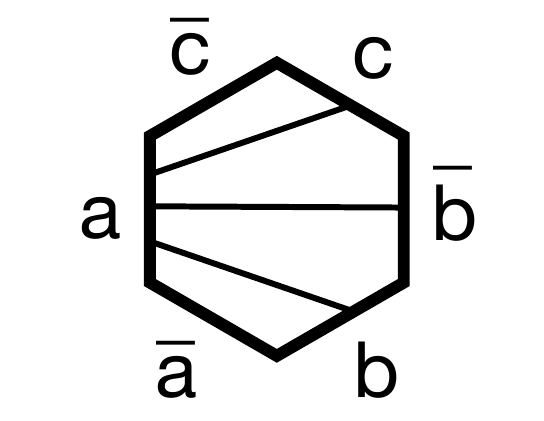}}

\bold{Proof of Theorem \ref{thm:quad}.}
Let $\cA = \{a,b,c,\abar,\bbar,\cbar\} \subset G$ be the alphabet consisting of the generators of $G$
and their inverses.   We can also regard $\cA$ as the cyclically ordered set of vertices of $H$ adjacent to $e$.
Let $P$ be a hexagon whose sides are labeled by the elements of $\cA$ in the same cyclic order; see Figure \ref{fig:abc}.
A {\em chord} $[u,v]$ is an unordered pair of distinct elements $u,v \mem \cA$
labeling non-adjacent sides of $P$.  For example, the chords incident to side $a$ are given by
$[a,b]$, $[a,\bbar]$ and $[a,c]$.  

The set of all chords $K$ consists of 9 elements.  We will also regard the 
$K$ as a basis for the vector space $\reals^K$.
Define a symmetric bilinear form on $\reals^K$ by $Q(k,k') = 1$ if the chords $k$ and $k'$ cross in $P$,
and $Q(k,k') = 0$ otherwise.  (By definition, two chords incident to the same side of $P$ do not cross.)
It is easily checked that $Q$ is an indefinite form of signature $(6,3)$,
by diagonalizing the given matrix.

Let $S = \{a,b,c\}$, and let
\begin{displaymath}
	\del : \reals^K \arrow \reals^S
\end{displaymath}
be the linear boundary map that
counts $+1$ when a generator occurs as an endpoint of a chord, and $-1$ when its inverse occurs.
(For example, $\del([a,b]) = a+b$, while $\del([a,\bbar]) = a-b$.)
Let
\begin{displaymath}
	\lambda : \reals^K \arrow \reals
\end{displaymath}
be the length function defined by summing the coordinates; that is, by setting $\lambda(k)=1$ for each chord $k$.
One can now readily verify that $Q|\Ker(\del)$ is positive--definite, and that for
any vector $v \mem \Ker(\del)$, we have
\begin{equation}
\label{eq:Q}
	Q(v,v) \ge \lambda(v)^2/3 .
\end{equation}
(The minimum is achieved on vectors that weight zero to all chords
that join opposite sides, and equal weights to the rest.)

Now let $w = g_1 \cdots g_n \mem G$ be a primitive, cyclically reduced word such that no pair
of adjacent letters $(g_{i-1},g_i)$ are equal.  (Here $g_n$ is adjacent to $g_1$.)
This condition, plus the fact that $w$ is reduced, insures that the sides of $P$
labeled by $\gbar_{i-1}$ and $g_i$ are not adjacent.
Thus we may associate to $w$ the sequence of chords 
\begin{displaymath}
	k_i = [\gbar_{i-1},g_i],
\end{displaymath}
$i=1,\ldots,n$.   

Let $v = \sum_1^n k_i \mem \reals^K$.  Then clearly $\delbar(v) = 0$ and $\lambda(v) = L(w) = n$.
Using formula (\ref{eq:I1}), it is also easy to see that
\begin{displaymath}
	I(w) \ge \frac{1}{2} \sum_{1 \le i,j \le n} Q(k_i,k_j) = \frac{1}{2} Q(v,v) .
\end{displaymath}
Indeed, when $Q(k_i,k_j) = 1$ the geodesics $\delta_i$ and $\delta_j$ cross transversally,
with degree $4$, at the identity element $e \mem H$, and hence we obtain a contribution
of $1/2$ to the sum in formula (\ref{eq:I1}).  Referring to inequality (\ref{eq:Q}), we obtain
the desired inequality
\begin{displaymath}
	I(w) \ge \lambda(v)^2/6 = L(w)^2 / 6.
\end{displaymath}
\qed

\bibliographystyle{math}
\bibliography{math}

\end{document}